# STALLINGS FOLDINGS AND SUBGROUPS OF FREE GROUPS

## ILYA KAPOVICH AND ALEXEI MYASNIKOV

ABSTRACT. We re-cast in a more combinatorial and computational form the topological approach of J.Stallings to the study of subgroups of free groups.

## 1. INTRODUCTION

The subgroup structure of free groups is a classical subject which goes back to the origins of Group Theory. The original approach developed by Nielsen was to treat this topic combinatorially, mainly using the technique of Nielsen transformations. Indeed, even to this day this method remains among the most powerful ones for working with subgroups of free groups. The development of algebraic topology and covering space theory in the 1940s suggested a different, more geometric approach. A free group $F$ can be identified with the fundamental group of a topological graph (which we may think of as a 1-complex). Then any subgroup of $F$ corresponds to a covering map from another graph to the original graph. The topological viewpoint was studied in detail by J.Stallings in a seminal paper [43]. In this work J.Stallings introduced an extremely useful notion of a *folding* of graphs. The ideas of Stallings have found many interesting applications (see, for example, [15], [16], [18], [44],,[45], [3], [41], [11], [13], [14], [48], [49], [50]). However, most of these applications, as far as the theory of free groups is concerned, were for the study of automorphisms of free groups. Thus, for example, the notion of a folding plays an important role in the construction of train-track maps due to M.Bestvina and M.Handel [6], [7].

In this paper we want to recast the ideas of Stallings' work [43] in a more combinatorial guise and apply them more systematically to the subgroup structure of free groups.

Namely, to every subgroup $H$ of a free group $F(X)$ we will associate a directed graph $\Gamma(H)$ whose edges are labeled by the elements of $X$. It turns out that this graph $\Gamma(H)$ (which is finite and very easy to construct if $H$ is finitely generated) carries all the essential information about the subgroup $H$ itself. Geometrically, the graph $\Gamma(H)$ represents the topological core of the covering space, corresponding to $H$, of the wedge of $\#X$ circles. Algebraically, $\Gamma(H)$ can also be viewed as the "essential part" of the relative coset Cayley graph of $G/H$ with respect to $X$. This last approach was used by C.Sims in [42]. Yet another view of $\Gamma(H)$ comes from Bass-Serre theory. Namely, $F(X)$ acts on its standard Cayley graph, which is a regular tree. There is a unique $H$-invariant subtree $T(H)$ which is minimal among $H$-invariant subtrees containing 1. Our $\Gamma(H)$ can be identified with the quotient $T(H)/H$.

However, in the present paper we will consider and study $\Gamma(H)$ almost exclusively as a combinatorial object. Namely, we will view it as a labeled directed graph and, also, as an automaton over the alphabet $X \cup X^{-1}$. There are several reasons why we believe this approach is interesting and worthwhile.

First and foremost, the graph-automaton $\Gamma(H)$ turns out to be an extremely useful and natural object when one considers various algorithmic and computational problems for free groups and their subgroups (e.g. the membership problem). In fact, most classical algorithms dealing with such problems usually use Nielsen methods and involve enumerating all elements of $F(X)$ of bounded length. Since $F(X)$ has exponential growth, all such algorithms require at least an exponential amount of time. However, it turns out that looking at $\Gamma(H)$ allows one to solve most of these problems much faster, ordinarily in polynomial time. In fact, much of this work is motivated by the desire to describe various algorithms implemented in the Computational Group Theory Software Package MAGNUS. It is also important that

---

the graph-automaton $\Gamma(H)$ is a canonical object which gives us some advantages as compared to the standard situation in the theory of automata and formal languages.

Second, we thought it worth-while to give an account of Stallings' ideas in a manner accessible to people with little topological or even group-theoretic background. There have been a number of papers where these ideas were applied to some specific subgroup properties of free groups. However, we wanted to give a self-contained, elementary and comprehensive treatment of the subject. For this reason we included complete and independent proofs of most basic facts, a substantial number of explicit examples and a wide assortment of possible applications.

Thus we re-prove many classical, well-known or folklore results about the subgroup structure of free groups. For instance, we prove the Marshal Hall Theorem and the Greenberg-Stallings Theorem as well as establish some well-known facts regarding ascending and descending chains of subgroups in free groups. In some instances we choose not to prove certain statements in the strongest possible form, as our real goal is to demonstrate the usefulness of the method.

However, there are quite a few results that appear to be new, or at least, have never been published. Thus in Section 11 we develop a mini-theory of "algebraic" and "free" extensions in the context of free groups. This technique is then applied in Section 13 to study malnormal closures and isolators of finitely generated subgroups of free groups.

This paper is based, in large measure, on the notes of a Group Theory Research Seminar, that has been running at the CUNY Graduate Center in the Spring and Fall semesters of 1999. Our special thanks go to Gilbert Baumslag, Toshiaki Jitsukawa, Gretchen Ostheimer, Gillian Elston, Robert Gilman, Vladimir Shpilrain, Sean Cleary, Lev Shneyerson, Fuh Ching-Fen, Dmitri Pechkin, Dmitri Bormotov, Alexei D. Myasnikov, Alexei Kvaschuk, Denis Serbin Katalin Bencsath, Arthur Sternberg and other regular participants of the seminar for the many lively, engaging and stimulating discussions.

## 2. Labeled graphs

**Definition 2.1** ($X$-digraph). Let $X = \{x_1, \ldots, x_N\}$ be a finite alphabet. By a $X$-*labeled directed graph* $\Gamma$ (also called a $X$-*digraph* or even just $X$-*graph*) we mean the following:

$\Gamma$ is a combinatorial graph where every edge $e$ has an arrow (direction) and is labeled by a letter from $X$, denoted $\mu(e)$. (Note that $\Gamma$ can be either finite or infinite.)

For each edge $e$ of $\Gamma$ we denote the origin of $e$ by $o(e)$ and the terminus of $e$ by $t(e)$. If $o(e) = t(e)$ then $e$ is a *loop*.

There is an obvious notion of a *morphism* between two $X$-digraphs. Namely, if $\Gamma$ and $\Delta$ are $X$-digraphs then a map $\pi : \Gamma \longrightarrow \Delta$ is called a *morphism* of $X$-digraphs, if $\pi$ takes vertices to vertices, directed edges to directed edges, preserves labels of directed edges and has the property that $o(\pi(e)) = \pi(o(e)), t(\pi(e)) = \pi(t(e))$ for any edge $e$ of $\Gamma$.

Examples of $X$-graphs, with $X = \{a, b, c\}$ are shown in Figure 1.

For the remainder of the paper, unless specified otherwise, $X$ will stand for a finite alphabet $X = \{x_1, \ldots, x_N\}$. We shall also denote $\Sigma = X \cup X^{-1}$.

**Convention 2.2.** Given a $X$-digraph $\Gamma$, we can make $\Gamma$ into an oriented graph labeled by the alphabet $\Sigma = X \cup X^{-1}$. Namely, for each edge $e$ of $\Gamma$ we introduce a formal inverse $e^{-1}$ of $e$ with label $\mu(e)^{-1}$ and the endpoints defined as $o(e^{-1}) = t(e)$, $t(e^{-1}) = o(e)$. The arrow on $e^{-1}$ points from the terminus of $e$ to the origin of $e$. For the new edges $e^{-1}$ we set $(e^{-1})^{-1} = e$.

The new graph, endowed with this additional structure, will be denoted by $\hat{\Gamma}$. In fact in many instances we will abuse notation by disregarding the difference between $\Gamma$ and $\hat{\Gamma}$.

The edge-set of $\hat{\Gamma}$ is naturally partitioned as $E\hat{\Gamma} = E\Gamma \cup \bar{E}\Gamma$. We will say that the edges $e$ of $\Gamma$ are *positively oriented* or *positive* in $\hat{\Gamma}$ and that their formal inverses $e^{-1}$ are *negatively oriented* or *negative* in $\hat{\Gamma}$. We will also denote $E^+\hat{\Gamma} = E\Gamma$ and $E^-\hat{\Gamma} = (E\Gamma)^{-1}$.

The use of $\hat{\Gamma}$ allows us to define the notion of a *path* in $\Gamma$. Namely, a *path* $p$ in $\Gamma$ is a sequence of edges $p = e_1, \ldots, e_k$ where each $e_i$ is an edge of $\hat{\Gamma}$ and the origin of each $e_i$ (for $i > 1$) is the terminus of $e_{i-1}$. In this situation we will say that the *origin* $o(p)$ of $p$ is $o(e_1)$ and the *terminus* $t(p)$ is $t(e_k)$. The *length* $|p|$ of this path is set to be $k$. Also, such a path $p$ has a naturally defined label $\mu(p) = \mu(e_1) \ldots \mu(e_k)$. Thus $\mu(p)$ is a word in the alphabet $\Sigma = X \cup X^{-1}$. Note that it is possible that $\mu(p)$ contains subwords of the form $aa^{-1}$ or $a^{-1}a$ for some $a \in X$.

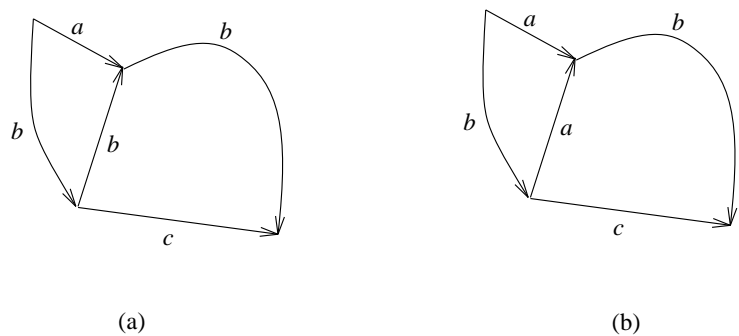

FIGURE 1. Folded and non-folded graphs. The labeling alphabet is $X = \{a, b, c\}$.

Also, if $v$ is a vertex of $\Gamma$, we will consider the sequence $p = v$ to be a path with $o(p) = t(p) = v$, $|p| = 0$ and $\mu(p) = 1$ (the empty word).

**Definition 2.3** (Folded graphs). Let $\Gamma$ be an $X$-digraph. We say that $\Gamma$ is *folded* if for each vertex $v$ of $\Gamma$ and each letter $a \in X$ there is at most one edge in $\Gamma$ with origin $v$ and label $a$ and there is at most one edge with terminus $v$ and label $a$.

The graph shown in Figure 1(a) is folded and the graph shown in Figure 1(b) is not folded. Note that $\Gamma$ is folded if and only if for each vertex $v$ of $\Gamma$ and each $x \in \Sigma = X \cup X^{-1}$ there is at most one edge in $\hat{\Gamma}$ with origin $v$ and label $x$. Thus in a folded $X$-digraph the degree of each vertex is at most $2\#(X)$.

The following notion of graph folding plays a fundamental role in this paper. Suppose $\Gamma$ is an $X$-digraph and $e_1, e_2$ are edges of $\Gamma$ with common origin and the same label $x \in \Sigma$. Then, informally speaking, *folding* $\Gamma$ at $e_1, e_2$ means identifying $e_1$ and $e_2$ into a single new edge labeled $x$. The resulting graph carries a natural structure of an $X$-digraph. A more precise definition is given below.

**Definition 2.4** (Folding of graphs). Let $\Gamma$ be an $X$-digraph. Suppose that $v_0$ is a vertex of $\Gamma$ and $f_1, f_2$ are two distinct edges of $\hat{\Gamma}$ with origin $v_0$ and such that $\mu(f_1) = \mu(f_2) = x \in \Sigma = X \cup X^{-1}$ (so that $\Gamma$ is not folded). Let $h_i$ be the positive edge of $\Gamma$ corresponding to $f_i$ (that is $h_i = f_i$ if $f_i$ is positive and $h_i = f_i^{-1}$ if $f_i$ is negative). Note that, depending on whether $x \in X$ or $x \in X^{-1}$, the edges $f_1$ and $f_2$ are either both positive or both negative in $\hat{\Gamma}$.

Let $\Delta$ be an $X$-digraph defined as follows.

The vertex set of $\Delta$ is the vertex set of $\Gamma$ with $t(f_1)$ and $t(f_2)$ removed and a new vertex $t_f$ added (we think of the vertices $t(f_1)$ and $t(f_2)$ as being identified to produce vertex $t_f$):

$$V\Delta = \big(V\Gamma - \{t(f_1), t(f_2)\}\big) \cup \{t_f\}$$

The edge set of $\Delta$ is the edge set of $\Gamma$ with the edges $h_1, h_2$ removed and a new edge $h$ added (we think of the edges $h_1$ and $h_2$ as being identified or "folded" to produce a new edge $h$):

$$E\Delta = (E\Gamma - \{h_1, h_2\}) \cup \{h\}$$

The endpoints and arrows for the edges of $\Delta$ are defined in a natural way. Namely, if $e \in E\Delta$ and $e \neq h$ (that is $e \in E\Gamma$, $e \neq h_i$) then

1. we put $o_\Delta(e) = o_\Gamma(e)$ if $o_\Gamma(e) \neq t(f_i)$ and $o_\Delta(e) = t_f$ if $o_\Gamma(e) = t(f_i)$ for some $i$;
2. we put $t_\Delta(e) = t_\Gamma(e)$ if $t_\Gamma(e) \neq t(f_i)$ and $t_\Delta(e) = t_f$ if $t_\Gamma(e) = t(f_i)$ for some $i$

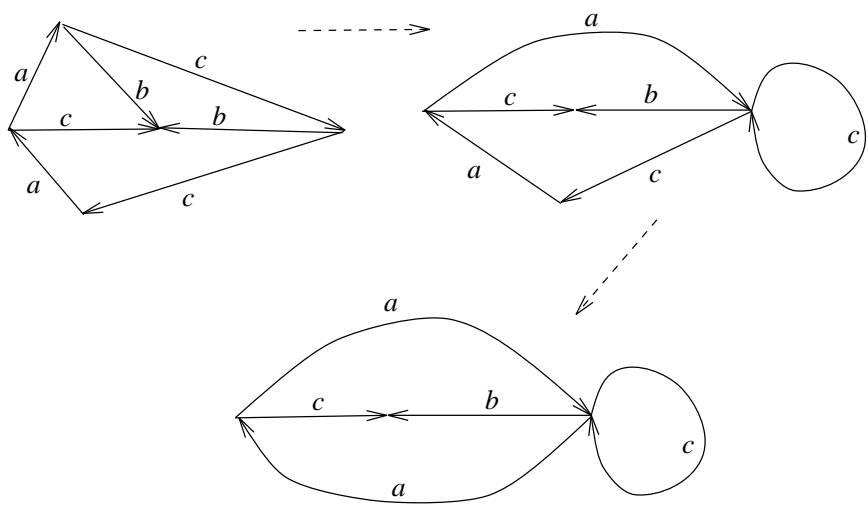

FIGURE 2. Folding of labeled graphs

For the edge $h$ we put $o_\Delta(h) = o_\Gamma(h)$ if $o_\Gamma(h) \neq t(f_1), t(f_2)$ and $o_\Delta(h) = t_f$ if $o_\Gamma(h) = t(f_i)$ for some $i$. Similarly, we put $t_\Delta(h) = t_\Gamma(h)$ if $t_\Gamma(h) \neq t(f_1), t(f_2)$ and $t_\Delta(h) = t_f$ if $t_\Gamma(h) = t(f_i)$ for some $i$.

We define labels on the edges of $\Delta$ as follows: $\mu_\Delta(e) = \mu_\Gamma(e)$ if $e \neq h$ and $\mu_\Delta(h) = \mu_\Gamma(h_1) = \mu_\Gamma(h_2)$.

Thus $\Delta$ is an $X$-digraph. In this situation we say that $\Delta$ is obtained from $\Gamma$ by a *folding* (or by *folding the edges $f_1$ and $f_2$*).

We suggest that the reader investigates carefully what happens in the above definition if one or both of the edges $f_1, f_2$ are loops. The notion of folding is illustrated in Figure 2.

We also need the following obvious statement summarizing some basic properties of folding.

**Lemma 2.5.** *Let $\Gamma_1$ be an $X$-digraph obtained by a folding from a graph $\Gamma$. Let $v$ be a vertex of $\Gamma$ and let $v_1$ be the corresponding vertex of $\Gamma_1$.*

*Then the following holds:*

1. *If $\Gamma$ is connected then $\Gamma_1$ is connected.*
2. *Let $p$ be a path from $v$ to $v$ in $\Gamma$ with label $w$. Then the edge-wise image of $p$ in $\Gamma_1$ is a path from $v_1$ to $v_1$ with label $w$.*
3. *Note that if an $X$-digraph $\Gamma$ is finite, then a folding always decreases the number of edges in $\Gamma$ by one.*

**Convention 2.6** (Reduced words and reduced paths)**.** Recall that a word $w$ in the alphabet $\Sigma = X \cup X^{-1}$ is said to be *freely reduced* if it does not contain a subword of the form $aa^{-1}$ or $a^{-1}a$ for $a \in X$.

A path $p$ in an $X$-digraph $\Gamma$ is said to be *reduced* if $p$ does not contain subpaths of the form $e, e^{-1}$ for $e \in E\dot{\Gamma}$.

For an alphabet $B$ we denote by $B^*$ the set of all words (including the empty word 1) in the alphabet $B$.

**Definition 2.7** (Language recognized by a digraph)**.** Let $\Gamma$ be an $X$-digraph and let $v$ be a vertex of $\Gamma$. We define *the language of $\Gamma$ with respect to $v$* to be:

$$L(\Gamma, v) = \{\mu(p) \,|\, p \text{ is a reduced path in } \Gamma \text{ from } v \text{ to } v\}$$

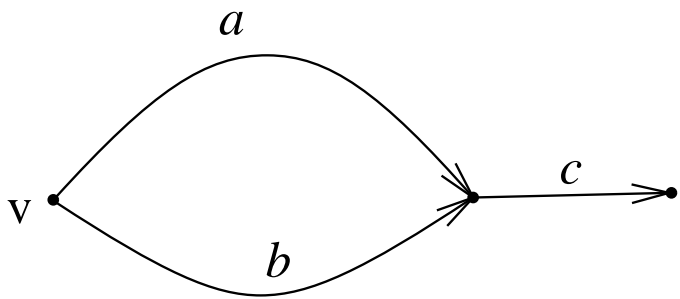

FIGURE 3. Language of a digraph. In this case the language $L(\Gamma, v)$ consists of all words of the form $(ab^{-1})^n$ where $n$ is any integer (possibly negative). Note that the edge labeled $c$ is irrelevant since there are no reduced paths from $v$ to $v$ passing through this edge.

If a word $w$ belongs to $L(\Gamma, v)$, we will also sometimes say that $w$ is *accepted* by $(\Gamma, v)$ (or just by $\Gamma$ if the choice of $v$ is fixed).

Thus $L(\Gamma, v)$ is a subset of $\Sigma^*$ where $\Sigma = X \cup X^{-1}$. An example of the language corresponding to a digraph is shown in Figure 3

Note that the words in $L(\Gamma, v)$ are not necessarily freely reduced. However, if the graph $\Gamma$ is folded, then all the words in $L(\Gamma, v)$ are freely reduced.

**Remark 2.8** (Labeled graphs versus standard automata)**.** It is possible and often useful to view the pair $(\hat{\Gamma}, v)$ as an automaton $M$ (possibly infinite) over the alphabet $\Sigma = X \cup X^{-1}$. The vertices of $\Gamma$ are states and the edges represent transitions. There is only one initial state, namely $v$ and only one accept state, also $v$. Then $\Gamma$ is folded if and only if the corresponding automaton $M$ is deterministic. It should be noted, that our definition of the language $L(\Gamma, v)$ is slightly different from the standard notion of the language $L(M)$ recognized by $M$. Namely, using the standard conventions of the automata theory, the language of $L(M)$ would consist of the labels of all (not just reduced) paths from $v$ to $v$ in $\Gamma$. However, the set of all freely reduced words in $\Sigma$ is easily seen to be a regular language (i.e. it is recognized by a finite state automaton). Since intersections of regular languages are regular, this implies that for a finite and folded $X$-digraph $\Gamma$, the language $L(\Gamma, v)$ is regular.

We formulate this obvious but important statement as a lemma.

**Lemma 2.9.** *Let $\Gamma$ be a folded $X$-digraph. Let $v$ be a vertex of $\Gamma$. Then all the words in the language $L(\Gamma, v)$ are freely reduced.*

## 3. Labeled graphs and subgroups of free groups

In this section we will explain how to every finitely generated subgroup of a free group one can associate a (reasonably) canonical labeled graph.

Let $X$ be a finite alphabet and let $\Sigma = X \cup X^{-1}$. If $w$ is a word in $\Sigma$, we will denote by $\overline{w}$ the freely reduced word in $\Sigma$ obtained from $w$ by performing all possible (if any) free reductions.

Recall that a *free group on $X$*, denoted $F(X)$ is the collection of all freely reduced words in $\Sigma$ (including the empty word 1). The multiplication in $F(X)$ is defined as

$$f \cdot g := \overline{fg}, \quad \text{for any } f, g \in F(X)$$

(Although not quite trivial, it can be shown that thus defined, the multiplication makes $F(X)$ a group [32]). The number of elements in $X$ is called *the rank of $F(X)$* and $X$ is referred to as *a free basis of $F(X)$*.

For any element $g \in F(X)$ we will denote by $|g|_X$ the length of the unique freely reduced $\Sigma$-word, representing $g$. If $w$ is a word in $\Sigma$ (which may or may not be freely reduced), we will denote by $|w|$ the length of $w$.

It turns out that the languages of directed $X$-graphs correspond to subgroups of $F(X)$.

Namely, the following simple but important statement holds.

**Proposition 3.1.** *Let $\Gamma$ be an $X$-digraph and let $v$ be a base-vertex of $\Gamma$. Then the set*

$$\overline{L} = \{\overline{w} \,|\, w \in L(\Gamma, v)\}$$

*is a subgroup of $F(X)$.*

*Proof.* Indeed, let $v_1, v_2 \in \overline{L}$. Then there are reduced paths $p_1$ and $p_2$ from $v$ to $v$ in $\Gamma$ such that the label of $p_i$ is $w_i$ and $\overline{w_i} = v_i$.

The concatenation $p_1 p_2$ is a path in $\Gamma$ from $v$ to $v$ which may or may not be path-reduced. Let $p$ be the reduced path obtained from $p_1 p_2$ by making all possible path reductions. This means that the label $w = \mu(p)$ is obtained from the word $w_1 w_2$ by performing several free reductions (even though $w = \mu(p)$ may not be freely reduced itself). Therefore $\overline{w} = \overline{w_1 w_2} = v_1 \cdot v_2 \in F(X)$. On the other hand $w$ is the label of a reduced path from $v$ to $v$ and therefore $w \in L(\Gamma, v) = L$ (by definition). Thus $v_1 \cdot v_2 \in \overline{L}$ and $\overline{L}$ is closed under multiplication.

It is easy to see that the inverse path $(p_1)^{-1}$ of $p_1$ has label $w_1^{-1}$. This implies that $\overline{L}$ is closed under taking inverses. Also, obviously $1 \in \overline{L}$.

Thus $\overline{L}$ is a subgroup of $F(X)$, as required. ∎

**Lemma 3.2.** *Suppose $\Gamma$ is a folded $X$-digraph. Then $L(\Gamma, v) = L = \overline{L}$ is a subgroup of $F(X)$.*

*Proof.* By Lemma 2.9 the labels of reduced paths in $\Gamma$ are already freely reduced. Therefore by Proposition 3.1 $L(\Gamma, v) = L = \overline{L}$ is a subgroup of $F(X)$. ∎

The following lemma shows that subgroups corresponding to finite labeled graphs are finitely generated.

**Lemma 3.3.** *Let $\Gamma$ be a connected $X$-digraph and let $v$ be a vertex of $\Gamma$. For each vertex $u \neq v$ of $\Gamma$ choose a reduced path $p_u$ in $\Gamma$ from $v$ to $u$. Put $p_v = v$, the path of length zero consisting just of vertex $v$. For each edge $e$ of $\hat{\Gamma}$ (whether positive or negative) put $p_e = p_{o(e)} e (p_{o(e)})^{-1}$ so that $p_e$ is a path in $\Gamma$ from $v$ to $v$. Denote $[e] = \overline{\mu(p_e)}$.*

*Then the subgroup $H = \overline{L(\Gamma, v)}$ of $F(X)$ is generated by the set*

$$S = \{[e] \,|\, \text{ where } e \text{ is a positive edge of } \Gamma\}$$

*In particular, if $\Gamma$ is finite, the subgroup $H$ is finitely generated.*

*Proof.* The path $p_e$ (if it is not already reduced) can be transformed into a reduced path $p'_e$ (from $v$ to $v$) by a series of path-reductions. Clearly $\mu(p'_e) \in L(\Gamma, v)$, $\overline{\mu(p'_e)} = \overline{\mu(p_e)}$ and hence $[e] = \overline{\mu(p_e)} \in \overline{L(\Gamma, v)} = H$. Thus $S \subseteq H$ and $\langle S \rangle \leq H$.

It remains to prove that any element of $H$ can be expressed as a product of elements of $S$ and their inverses. Note that by definition $p_{e^{-1}} = (p_e)^{-1}$ and so $[e^{-1}] = [e]^{-1}$. Thus it suffices to show that any element of $H$ can be expressed as a product of elements $[e]$, where $e$ is an edge of $\hat{\Gamma}$.

Let $h \in H$, $h \neq 1$. Then there is a nontrivial reduced path $p$ from $v$ to $v$ with label $w$ such that $\overline{w} = h$. Let $p = e_1, \ldots, e_k$, where $e_i$ are edges of $\hat{X}$. Further let $v_1 = v$, $v_{k+1} = v$ and let $v_i$ be the initial vertex of $e_i$ (and therefore the terminal vertex of $e_{i+1}$). Consider now the path

$$p' = p_{e_1} \ldots p_{e_k} = p_{v_1} e_1 p_{v_2}^{-1} p_{v_2} e_2 p_{v_3}^{-1} \ldots p_{v_k} e_k p_{v_{k+1}}^{-1}$$

It is obvious that the path $p'$ can be transformed by path reductions into the path $p$. Therefore $\overline{\mu(p')} = \overline{\mu(p)} = h$. On the other hand

$$\overline{\mu(p')} = \overline{p_{e_1}} \ldots \overline{p_{e_k}} = [e_1] \ldots [e_k] \in \langle S \rangle$$

Thus $H = \langle S \rangle$ as required. ∎

The following simple but important observation plays a crucial role in this paper.

**Lemma 3.4.** *Let $\Gamma$ be an $X$-digraph and let $\Gamma'$ be an $X$-digraph obtained from $\Gamma$ by a single folding. Further let $v$ be a vertex of $\Gamma$ and let $v'$ be the corresponding vertex of $\Gamma'$.*

*Then $\overline{L(\Gamma, v)} = \overline{L(\Gamma', v')}$.*

*Proof.* Suppose $\Gamma'$ is obtained from $\Gamma$ by folding two edges $e_1, e_2$ in $\hat{\Gamma}$ which have the same initial vertex $u$ and the same label $x \in X$. The edges $e_1, e_2$ are folded into an edge $e$ of $\Gamma'$ labeled $x$ and with origin $u'$.

Suppose $p$ is a reduced path in $\Gamma$ from $v$ to $v$, so that $\mu(p) \in L(\Gamma, v)$. The image of $p$ in $\Gamma'$ is a path $p'$ from $v'$ to $v'$ with the same label as that of $p$, that is $\mu(p) = \mu(p')$. However, $p'$ need not be path-reduced. Namely, $p'$ is path-reduced if and only if $p$ does not contain any subpaths of the form $e_2^{-1}, e_1$ or $e_1^{-1}, e_2$. Let $p''$ be the path obtained from $p'$ by performing all possible path reductions in $\Gamma'$. Then $\overline{\mu(p)} = \overline{\mu(p')} = \overline{\mu(p'')}$ and $\mu(p'') \in L(\Gamma, v')$. Thus we have shown that $\overline{L(\Gamma, v)} \subseteq \overline{L(\Gamma', v')}$.

Suppose now that $p'$ is an arbitrary reduced path in $\Gamma'$ from $v'$ to $v'$. We claim that there is a reduced path in $\Gamma$ from $v$ to $v$ with exactly the same label as that of $p'$. We will construct this path explicitly.

The occurrences of $e^{\pm 1}$ (if any) subdivide $p'$ into a concatenation of the form:

$$p' = p_0 f_0 p_1 f_1 \ldots f_k p_{k+1}$$

where $f_i = e^{\pm 1}$ and the paths $p_i$ do not involve $e^{\pm 1}$.

Suppose that for some $i$ we have $f_i = e$. Since $p_i$ and $p_{i+1}$ do not involve the edge $e$, they can also be considered as paths in $\Gamma$. Moreover, by the definition of folding, in the graph $\Gamma$ the terminal vertex of $p_i$ is joined with the initial vertex of $p_{i+1}$ by either the edge $e_1$ or the edge $e_2$. We denote this edge by $d_i$ (so that $d_i \in \{e_1, e_2\}$). Note that now $p_i d_i p_{i+1}$ is a reduced path in $\Gamma$ with the same label as the path $p_i f_i p_{i+1}$ in $\Gamma'$.

Similarly, if for some $i$ we have $f_i = e^{-1}$, we can find $d_i \in \{e_1^{-1}, e_2^{-1}\}$ such that $p_i d_i p_{i+1}$ is a reduced path in $\Gamma$ with the same label as the path $p_i f_i p_{i+1}$ in $\Gamma'$.

Then

$$p = p_0 d_0 p_1 \ldots d_k p_{k+1}$$

is a reduced path in $\Gamma$ from $v$ to $v$ with the same label as $p'$. Thus $\mu(p') \in L(\Gamma, v)$, $\overline{\mu(p')} \in \overline{L(\Gamma, v)}$ and therefore $\overline{L(\Gamma', v')} \subseteq \overline{L(\Gamma, v)}$.

Hence $\overline{L(\Gamma', v')} \subseteq \overline{L(\Gamma, v)}$ and the lemma is proved. $\qquad\square$

Before investigating further group-theoretic properties of labeled graph we need to introduce the following important geometric concept.

**Definition 3.5** (Core graphs). Let $\Gamma$ be an $X$-digraph and let $v$ be a vertex of $\Gamma$. Then the *core of $\Gamma$ at $v$* is defined as:

$$Core(\Gamma, v) = \cup \{p \mid \text{ where } p \text{ is a reduced path in } \Gamma \text{ from } v \text{ to } v \}$$

It is easy to see that $Core(\Gamma, v)$ is a connected subgraph of $\Gamma$ containing $v$. If $Core(\Gamma, v) = \Gamma$ we say that $\Gamma$ is a *core graph with respect to $v$*.

**Example 3.6.** The graph shown in Figure 3 is not a core graph with respect to $v$. However, it is a core graph with respect to the terminus of the edge labeled $c$. Both graphs shown in Figure 1 are core graphs with respect to any of their vertices.

The following lemma lists some obvious properties of core graphs.

**Lemma 3.7.** *Let $\Gamma' = Core(\Gamma, v)$. Then:*

1. *The subgraph $\Gamma'$ of $\Gamma$ is connected and contains the vertex $v$.*
2. *The graph $\Gamma'$ has no degree-one vertices, except possibly for the vertex $v$.*
3. *The languages of $\Gamma$ and $\Gamma'$ at $v$ coincide, that is $L(\Gamma, v) = L(\Gamma', v)$.*

Since languages of folded graphs consist of freely reduced words, Lemma 3.4 implies that for a folded $X$-digraph $\Gamma$ the language $L(\Gamma, v)$ is a subgroup of $F(X)$. In fact, by the properties of core graphs, if $\Gamma' = Core(\Gamma, v)$ then $L(\Gamma, v) = L(\Gamma', v)$ is a subgroup of $F(X)$. We will now show that any finitely generated subgroup of $F(X)$ (and, later, even infinitely generated) can be obtained in such a way.

**Proposition 3.8.** *Let $H$ be a finitely generated subgroup of $F(X)$. Then there exists a finite $X$-digraph $\Gamma$ and a vertex $v$ of $\Gamma$ such that:*

1. *the graph $\Gamma$ is folded and connected;*
2. *all vertices of $\Gamma$, except possibly for vertex $v$, have degree greater than one;*
3. *the degree of each vertex in $\Gamma$ is at most $2\#(X)$;*
4. *the graph $\Gamma$ is a core graph with respect to $v$;*
5. *the language of $\Gamma$ is equal to $H$, that is $L(\Gamma, v) = H$.*

*Proof.* If $H = 1$ then the graph consisting of a single vertex obviously satisfies all the requirements of the proposition. Assume now that $H$ is a nontrivial finitely generated subgroup of $F(X)$.

Let $H$ be generated by the elements $h_1, \ldots, h_m$ (where we think of each $h_i$ as a freely reduced word in $\Sigma = X \cup X^{-1}$).

We define an $X$-digraph $\Gamma_1$ as follows. The graph $\Gamma_1$ is a wedge of $m$ circles wedged at a vertex called $v_1$. The $i$-th circle is subdivided into $|h_i|$ edges which are oriented and labeled by $X$ so that the label of the $i$-th circle (as read from $v_1$ to $v_1$) is precisely the word $h_i$.

Then any freely reduced word in $h_1, \ldots, h_m$ is the label of a reduced path in $\Gamma_1$ from $v$ to $v$. The converse is also obviously true. Thus $\overline{L(\Gamma_1, v_1)} = \langle h_1, \ldots, h_m \rangle = H$. Note also that $\Gamma_1$ is connected by construction and has no vertices of degree one.

We define a sequence of graphs $\Gamma_1, \Gamma_2, \ldots$ inductively as follows (see Figure 4 for a specific example). Suppose $\Gamma_i$ is already constructed. It is folded, terminate the sequence. Otherwise, let $\Gamma_{i+1}$ be obtained from $\Gamma_i$ by a folding.

Since $\Gamma_1$ is finite and a folding decreases the number of edges, this sequence terminates in finitely many steps with a folded graph $\Gamma_k$. Moreover, since a folding of a connected graph is connected, the graphs $\underline{\Gamma_1, \ldots, \Gamma_k}$ are connected. Let $v_i$ be the image of $v_1$ in the graph $\Gamma_i$. By Lemma 3.4 we have $H = \overline{L(\Gamma_1, v_1)} = \overline{L(\Gamma_k, v_k)} = L(\Gamma_k, v_k)$, the last equality implied by the fact that $\Gamma_k$ is folded.

Put $\Gamma = \Gamma_k$ and $v = v_k$. We have already shown that $L(\Gamma, v) = H$ and that $\Gamma$ is a connected folded finite graph. Since $\Gamma$ is folded, the degrees of its vertices are at most $2\#(X) = \#(X)$. We claim that $\Gamma$ is in fact a core graph with respect to $v$. Indeed, suppose this is not so and there is a degree-one vertex $u$ of $\Gamma$ which is different from $v$. Let $e$ be the unique edge of $\hat{\Gamma}$ with terminus $u$ and let $x \in \Sigma$ be the label of $e$. There exists an edge $e_1$ of the graph $\Gamma_1$ such that the image of $e_1$ in $\Gamma = \Gamma_k$ is $e$ (so that the label of $e_1$ is also $x$. Let $u_1$ be the terminus of $e_1$. Since $u \neq v = v_k$, we have $u_1 \neq v_1$. Recall the explicit construction of $\Gamma_1$ as the wedge of circles labeled $h_1, \ldots, h_m$. It follows from this construction that there is a path $p_1$ in $\Gamma_1$ from $v_1$ to $v_1$ which passes through the edge $e_1$ and has freely reduced label $w$ (namely, we can take $w = h_i^{\pm 1}$ for some $i$). Let $p$ be the image of the path $p_1$ in $\Gamma_k$ after performing all $k - 1$ foldings. Then $p$ is a path from $v$ to $v$ passing through $e$ and with a freely reduced label $w$. However the vertex $u$ (different from $v$) is of degree one. Therefore any path from $v$ to $v$ passing through $e$ contains a subpath $e, e^{-1}$ and hence cannot have a freely reduced label. This gives us a contradiction. Thus $Core(\Gamma, v) = \Gamma$ and the pair $(\Gamma, v)$ satisfies all the requirements of the proposition. $\square$

**Lemma 3.9.** *Let $\Gamma$ be a connected folded $X$-digraph which is a core graph with respect to some vertex $v$. Let $H = L(\Gamma, v)$. Then:*

1. *For any initial segment $w$ of a freely reduced word $h \in H$ there exists a unique reduced path $p$ in $\Gamma$ with origin $v$ and label $w$.*
2. *For any reduced path $p$ in $\Gamma$ with origin $v$ and label $w$, the word $w$ is a subword of some freely reduced word $h \in H$.*

*Proof.* (1) The uniqueness of $p$ follows from the fact that $\Gamma$ is folded. The existence of $p$ is also obvious. Indeed, we can take $p$ to be the initial segment of length $|w|$ of the path in $\Gamma$ from $v$ to $v$ with label $h$.

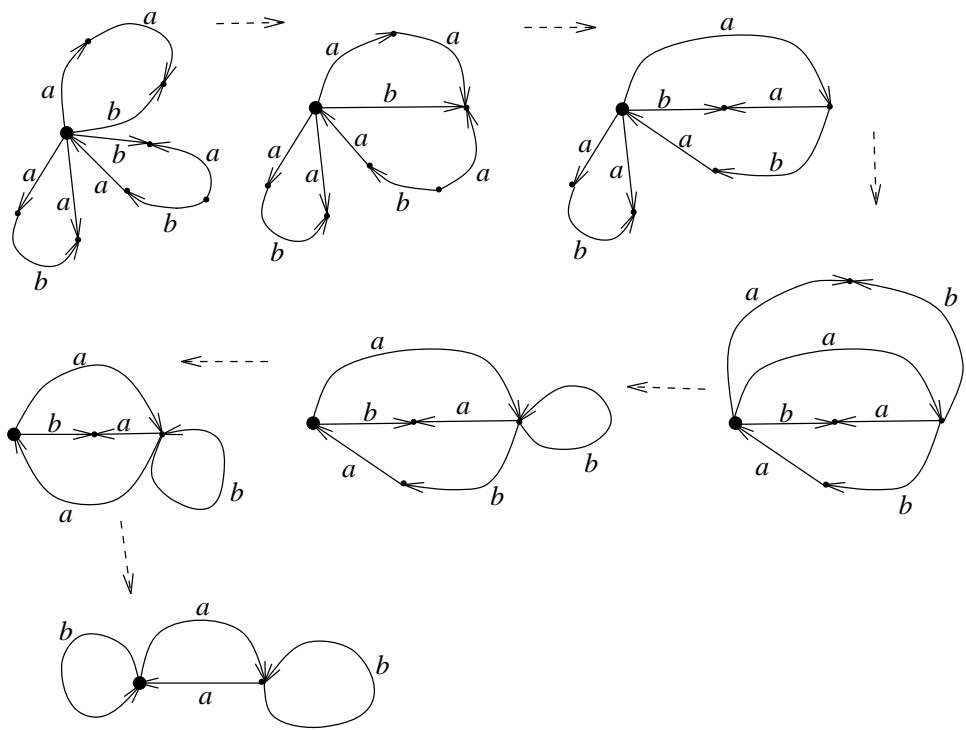

FIGURE 4. Constructing the graph associated to a subgroup. In this example we start with a subgroup $H = \langle a^2b, ba^{-1}ba, aba^{-1} \rangle \leq F(a, b)$. The highlighted vertices represent the images of the original base-vertex $v_1$ in the intermediate graphs

(2) Let $p$ be a path in $\Gamma$ with origin $v$ and label $w$. If the terminal vertex of $p$ is $v$ then $w \in H$ and the statement of the lemma is obvious. Suppose now that the terminal vertex $u$ of $p$ is different from $v$. Since $\Gamma$ is a core graph with respect to $v$, the degree of $u$ is at least two. Let $e$ be an edge with origin $u$ such that $e^{-1}$ is different from the last edge in $p$. Let $p'$ be a reduced path from $t(e)$ to $v$ in $\Gamma$. Then $\alpha = pep'$ is a reduced path from $v$ to $v$. Therefore $\mu(\alpha) = h \in H$. It is clear by construction that $w$ is an initial segment of $h$. □

## 4. Morphisms of labeled graphs

Recall that in Definition 2.1 we defined the notion of a morphism between two $X$-digraphs.
First we prove the following obvious but important lemma.

**Lemma 4.1** (Morphisms and subgroups)**.** *Let $\pi : \Gamma \longrightarrow \Gamma'$ be a morphism of $X$-digraphs such that $\pi(v) = v'$ for some vertex $v$ of $\Gamma$. Suppose that $\Gamma$ and $\Gamma'$ are folded. Put $K = L(\Gamma, v)$ and $H = L(\Gamma', v')$. Then $K \leq H$.*

*Proof.* Suppose $p$ is an arbitrary reduced path in $\Gamma$ from $v$ to $v$ and let $w$ be the label of $p$. Note that $w$ is freely reduced since $\Gamma$ is folded and that $w \in L(\Gamma, v) = K$.

Since $\pi$ is a morphism of $X$-digraphs, $\pi(p)$ is a path in $\Gamma'$ from $v'$ to $v'$ with label $w$. Since $w$ is freely reduced, the path $\pi(p)$ is path-reduced and therefore $w \in L(\Gamma', v') = H$. The path $p$ was chosen arbitrarily and therefore $K \leq H$. □

**Lemma 4.2** (Uniqueness of morphisms for folded graphs). *Let $\Gamma$ and $\Delta$ be connected $X$-digraphs such that $\Delta$ is folded. Let $v$ be a vertex of $\Gamma$ and let $u$ be a vertex of $\Delta$. Then there exists at most one morphism of $X$-digraphs $f : \Gamma \longrightarrow \Delta$ such that $f(v) = u$.*

*Proof.* Suppose $f, j : \Gamma \longrightarrow \Delta$ are such that $f(v) = j(v) = u$. Let $v'$ be an arbitrary vertex of $\Gamma$. Since $\Gamma$ is connected, there exists a reduced path $p$ in $\Gamma$ from $v$ to $v'$. Denote the label of $p$ by $w$. Since $\Delta$ is folded, there exists at most one path in $\Delta$ with origin $u$ and label $w$. On the other hand both $f(p)$ and $j(p)$ are such paths and they have terminal vertices $f(v')$ and $j(v')$ respectively. Therefore $f(v') = j(v')$. Since $v' \in V\Gamma$ was chosen arbitrarily, we have shown that $f$ and $j$ coincide on the vertex set of $\Gamma$. Since $\Gamma$ is folded, it easily follows that $f$ and $j$ coincide on the set of edges as well. $\qquad\square$

**Proposition 4.3.** *Let $F(X)$ be a free group with finite basis $H$. Let $K \leq H \leq F(X)$ be subgroups of $F(X)$. Suppose $(\Gamma_1, v_1)$ and $(\Gamma_2, v_2)$ are connected folded based $X$-digraphs such that $\Gamma_i$ is a core graph with respect to $v_i$ and $L(\Gamma_1, v_1) = K$ and $L(\Gamma_2, v_2) = H$.*

*Then there exists a unique morphism of $X$-digraphs $\pi : \Gamma_1 \longrightarrow \Gamma_2$ such that $\pi(v_1) = v_2$.*

*Proof.* The uniqueness of $\pi$ follows from Lemma 4.2. Thus it suffices to show that such $\pi$ exists.

We will construct $\pi$ explicitly as follows. Let $v$ be a arbitrary vertex of $\Gamma_1$. Choose any reduced path $p_v$ in $\Gamma$ from $v_1$ to $v$ and let $w$ be the label of $p_v$. By Lemma 3.9 since $L(\Gamma_1, v_1) = K$ and $\Gamma_1$ is a folded core graph with respect to $v_1$, the word $w$ is an initial segment of some freely reduced word belonging to $K$. Since $K \leq H$ and $H = L(\Gamma_2, v_2)$ Lemma 3.9 also implies that there is a unique path $q_v$ in $\Gamma_2$ with origin $v_2$ and label $w$. We set $\pi(v)$ to be the terminal vertex of $q_v$. We claim that this definition of $\pi(v)$ does not depend on the choice of $p_v$. Indeed, let $p'_v$ be another reduced path in $\Gamma_1$ from $v_1$ to $v$ and let $w'$ be the label of $p'_v$. Again, let $q'_v$ be the path in $\Gamma_2$ with origin $v_2$ and label $w'$. Denote the terminal vertex of $p_v$ by $u$ and the terminal vertex of $p'_v$ by $u'$. Note that $p_v(p'_v)^{-1}$ is a loop in $\Gamma_1$ at $v_1$ with label $w(w')^{-1}$. Hence $\overline{w(w')^{-1}} = k \in K \leq H$. Let $r$ be the reduced path in $\Gamma_1$ from $v_1$ to $v_1$ with label $k$. Then $rp'_v$ is a path in $\Gamma_1$ from $v_1$ to $u'$. On the other hand the path-reduced form of $rp'_v$ has label $k \cdot w' = w$. Since $\Gamma_2$ is folded, there is no more than one path in $\Gamma_2$ with origin $v_1$ and label $w$, namely the path $p_v$. Since path-reductions do not change end-points of paths, this means that $u' = v'$, as required.

We now extend $\pi$ to the edges of $\Gamma_1$. Let $e$ be an edge of $\Gamma_1$ with label $x$. Since $\Gamma_2$ is folded, there is at most one edge in $\Gamma_2$ with label $x$ and origin $\pi(o(e))$. We claim that such an edge in fact exists. Since $\Gamma_1$ is a core graph with respect to $1_K$, there is a reduced path $p = p'ep''$ from $v_1$ to $v_1$ passing through $e$. Let $z' = \mu(p'), z'' = \mu(p'')$ so that $z = \mu(p'') = z'xz''$. By assumption on $\Gamma_1$ we have $z \in K$. However, $K \leq H$, so that $z \in H$ as well. Hence there exists a unique path $q$ is $\Gamma_2$ with origin $v_2$ and label $z = z'xz''$. This $q$ has the form $q = y'e'y''$ where the label of $q'$ is $z'$, the label of edge $e'$ is $x$ and the label of $y''$ is $z''$. By our construction of $\pi$ the terminus of $y'$ is $\pi(o(e))$. Hence the edge $e'$ has label $x$ and origin $\pi(o(e))$ in $\Gamma_2$. We put $\pi(e) = e'$.

Thus we have established that $\pi$ exists. $\qquad\square$

**Lemma 4.4** (Morphism factorization). *Let $f : \Gamma \longrightarrow \Delta$ be an epimorphism of $X$-digraphs. Then $f$ can be decomposed as $f = f_2 \circ f_1$ where $f_1 : \Gamma \longrightarrow \Gamma'$ is injective on the edge set of $\Gamma$ and $f_2 : \Gamma' \longrightarrow \Delta$ is a bijection between the vertex sets of $\Gamma'$ and $\Delta$.*

*Proof.* Let $f : \Gamma \longrightarrow \Delta$ be an epimorphism of $X$-digraphs. Define a graph $\Gamma'$ as follows. We give an informal description of $\Gamma'$, $f_1$ and $f_2$ and leave the technicalities to the reader. The map $f_1$ consists in identifying those vertices (but not edges) of $\Gamma$ which are identified by $f$. Thus $\Gamma' = f_1(\Gamma)$ is obtained by collapsing some subsets of $V\Gamma$ into single vertices.

Note now that if $e, h$ are edges of $\Gamma$ with the same image in $\Delta$ under $f$, then the origins of $f_1(e), f_1(h)$ are the same and the termini of $f_1(e), f_1(h)$ are the same. The morphism $f_2$ consists in identifying all such $f_1(e), f_1(h)$. Thus $f_2$ is in fact an identity map on the set of vertices of $\Gamma'$. Moreover, $f_2$ could be considered as a "generalized folding". In fact, $f_2$ is a composition of several foldings if $\Gamma$ is finite. $\qquad\square$

**Convention 4.5.** A *based $X$-digraph* is an $X$-digraph $\Gamma$ with a marked vertex $v$, called the *base-vertex*. Such a based digraph is denoted $(\Gamma, v)$. A *morphism* of based $X$-digraphs $\pi : (\Gamma, v) \longrightarrow (\Delta, u)$ is a morphism of the underlying $X$-digraphs $\pi : \Gamma \longrightarrow \Delta$ such that $\pi(v) = u$. By Lemma 4.2 if two folded connected based digraphs $(\Gamma, v)$ and $(\Delta, u)$ are isomorphic then there exists a unique isomorphism $\pi : (\Gamma, v) \longrightarrow (\Delta, u)$. Therefore we can identify $(\Gamma, v)$ and $(\Delta, u)$ via $\pi$. In this situation we will sometimes write $(\Gamma, v) = (\Delta, u)$.

## 5. Definition of subgroup graphs

In Proposition 3.8 we saw that one can associate to any finitely generated subgroup $H \le F(X)$ a finite folded connected $X$-digraph $\Gamma$ which is a core graph with respect to some vertex $v$ and such that $L(\Gamma, v)$. In this section we show that the same can be done for an arbitrary subgroup $H \le F(X)$ (although the graph will not be finite if $H$ is not finitely generated) and the pair $(\Gamma, v)$ is unique for a fixed $H$.

**Theorem 5.1.** *Let $F(X)$ be a free group with finite basis $X$. Let $H \le F(X)$ be a subgroup of $F(X)$. Then there exist a connected folded $X$-digraph $\Gamma$ and a vertex $v$ of $\Gamma$ such that:*

1. *the graph $\Gamma$ is a core graph with respect to $v$;*
2. *$L(\Gamma, v) = H$.*

*Proof.* First we will construct an $X$-digraph $\Delta$ as follows.

The vertex set of $\Delta$ is the set of left cosets of $H$ in $F(X)$:

$$V\Delta = \{Hf \,|\, f \in F(X)\}.$$

For two cosets $Hf$ and $Hg$ and a letter $x \in X$ we introduce a directed edge with origin $Hf$, terminus $Hg$ and label $x$ whenever $Hfx = Hg$. This defines an $X$-digraph $\Delta$. Put $v = H$ to be the coset of the identity element $1 \in F(X)$.

Note that $\Delta$ is a connected graph. Indeed, suppose $f = x_1 \ldots x_n$ is an arbitrary nontrivial freely reduced word in $X^{\pm 1}$ (here $x_i \in X \cup X^{-1}$). Then for each $i = 1, \ldots, n-1$ there is an edge in $\hat{\Delta}$ from $Hx_1 \ldots x_i$ to $Hx_1 \ldots x_i x_{i+1}$ with label $x_{i+1}$. Therefore there is a path in $\hat{\Delta}$ from $v = H1$ to $Hf$ with label $f$, and so $\Delta$ is connected.

Observe now that $\Gamma$ is folded. Suppose this is not the case and there is a vertex $u = Hf$ with two distinct edges $e_1$, $e_2$ which have the same label $x \in X \cup X^{-1}$. Let $u_1$ and $u_2$ be the terminal vertices of $e_1$ and $e_2$ accordingly. By definition of $\Delta$ this means that $u_1 = fxH = u_2$. However, by the construction of $\Delta$, for any pair of vertices $v_1, v_2$ and any $x \in X \cup X^{-1}$ there is at most one edge from $v_1$ to $v_2$ with label $x$. Thus $e_1 = e_2$, contrary to our assumptions.

We claim that $L(\Delta, v) = H$. Indeed, suppose $f$ is a nontrivial freely reduced word with is a label of a loop $p$ at $v$ in $\Delta$. Again let $f = x_1 \ldots, x_n$ where $x_i \in X \cup X^{-1}$. The definition of $\Delta$ implies that for every $i = 1, \ldots, n$ the terminal vertex of the initial segment of $p$ with label $x_1 \ldots x - i$ is the coset $Hx_1 \ldots x_i$. In particular the terminal vertex of $p$ is $Hx_1 \ldots x_n = Hf$. On the other hand the terminal vertex of $p$ is $H$ by assumption. Therefore $Hf = H$ and so $f \in H$. Thus we have shown that $L(\Delta, v) \subseteq H$.

Suppose now that $f$ is an arbitrary nontrivial freely reduced word such that $f \in H$. As we have seen when proving the connectivity of $\Delta$, there exists a path $p$ in $\Delta$ with origin $v = H$ and label $f$. Obviously $p$ is reduced since $f$ is freely reduced and $\Delta$ is folded. Let $u$ be the terminal vertex of $p$. By the construction of $p$ we have $u = Hf = H = v$, since $f \in H$. Thus $p$ is a reduced loop at $v$ with label $f$. Hence $f \in L(\Delta, v)$. Since $f \in H$ was chosen arbitrarily, we have proved that $H \subseteq L(\Delta, v)$.

Thus $H \subseteq L(\Delta, v) \subseteq H$ and therefore $H = L(\Delta, v)$.

Now put $\Gamma = Core(\Delta, v)$. It is obvious that $\Gamma$ satisfies all the requirements of Theorem 5.1. $\qquad \square$

In a seminal paper [42] C.Sims uses a similar approach to that used in the proof of Theorem 5.1 to study subgroups of free groups via "important" parts of their coset graphs.

**Theorem 5.2.** *Let $F(X)$ be a free group with finite basis $X$ and let $H \le F(X)$ be a subgroup of $F(X)$. Suppose $(\Gamma_1, v_1)$ and $(\Gamma_2, v_2)$ are connected folded based $X$-digraphs such that $\Gamma_i$ is a core graph with respect to $v_i$ and $L(\Gamma_1, v_1) = H = L(\Gamma_2, v_2)$.*

*Then there exists a unique isomorphism of $X$-digraphs $\pi : \Gamma_1 \longrightarrow \Gamma_2$ such that $\pi(v_1) = v_2$.*

*Proof.* The uniqueness of $\pi$ follows from Lemma 4.2.

Since $H \leq H$, by Proposition 4.3 there is a morphism $\pi : \Gamma_1 \longrightarrow \Gamma_2$ such that $\pi(v_1) = v_2$. We claim that $\pi$ is an isomorphism of $X$-digraphs. Indeed, Proposition 4.3 implies that there is a morphism $\pi' : \Gamma_2 \longrightarrow \Gamma_1$ such that $\pi'(v_2) = v_1$. Therefore $\pi' \circ \pi : \Gamma_1 \longrightarrow \Gamma_1$ is a morphism such that $(\pi' \circ \pi)(v_1) = v_1$. By Lemma 4.2 there is at most one such morphism, namely the identity map on $\Gamma_1$. Thus $\pi' \circ \pi = Id_{\Gamma_1}$. By a symmetric argument $\pi \circ \pi' = Id_{\Gamma_2}$. Therefore $\pi$ is an isomorphism as required.

**Definition 5.3** (Subgroup graph). Let $H \leq F(X)$ be a subgroup of $F(X)$. Then by Theorem 5.1 and Theorem 5.2 there exists a unique (up to a canonical isomorphism of based $X$-digraphs) based $X$-digraph $(\Gamma, v)$ such that

1. the graph $\Gamma$ is folded and connected;
2. the graph $\Gamma$ is a core graph with respect to $v$;
3. the language of $\Gamma$ with respect to $v$ is $H$, that is $L(\Gamma, v) = H$.

In this situation we call $\Gamma$ the *subgroup graph* (or a *subgroup graph-automaton*) of $H$ with respect to $X$ and denote it by $\Gamma(H)$ (or $\Gamma^X(H)$). The base-vertex $v$ is denoted $1_H$.

**Lemma 5.4.** *Let $H \leq F(X)$ be a subgroup of $F(X)$. Then the graph $\Gamma(H)$ is finite if and only if $H$ is finitely generated.*

*Proof.* If $H$ is finitely generated then $\Gamma(H)$ is finite by Proposition 3.8 and Proposition 4.3. If $\Gamma(H)$ is finite then $H$ is finitely generated by Lemma 6.1. $\qquad \blacksquare$

## 6. Spanning trees and free bases of subgroups of free groups

Lemma 3.3 provides a generating set for the subgroup corresponding to a $X$-digraph. However, one would like to be able to find a free basis or even a free Nielsen reduced basis for this subgroup. We will show next that both these goals can be easily accomplished using the original $X$-digraph.

Recall that in a connected graph a subgraph is said to be *a spanning tree* if this subgraph is a tree and it contains all vertices of the original graph. If a graph $T$ is a tree then for any two vertices $u$, $u'$ of $T$ there is a unique reduced path in $T$ from $u$ to $u'$ which will be denoted by $[u, u']_T$.

**Lemma 6.1.** *Let $\Gamma$ be a folded $X$-digraph and let $v$ be a vertex of $\Gamma$. Let $T$ be a spanning tree of $\Gamma$. Let $T^+$ be the set of those positive edges of $\Gamma$ which lie outside of $T$. For each $e \in T^+$ put $p_e = [v, o(e)]_T e[t(e), v]$ (so that $p_e$ is a reduced path from $v$ to $v$ and its label is a freely reduced word in $\Sigma = X \cup X^{-1}$.) Also for each $e \in T^+$ put $[e] = \mu(p_e) = \overline{\mu(p_e)}$. Denote*

$$Y_T = \{[e] \mid e \in T^+\}.$$

*Then $Y_T$ is a free basis for the subgroup $H = L(\Gamma, v)$ of $F(X)$.*

*Proof.* We can extend the definition of $p_e$ and $[e]$ for all edges of $\hat{\Gamma}$ (whether positive or negative and whether inside or outside of $T$) by putting $p_e = [v, o(e)]_T e[t(e), v]$ and $[e] = \overline{\mu(p_e)}$. Note that in this case $p_{e^{-1}} = (p_e)^{-1}$ and $[e^{-1}] = [e]^{-1}$.

It is easy to see that if $e \in E^+ T$ then $p_e = [v, o(e)]_T e[t(e), v]$ is a path that can be transformed by path-reductions into the trivial path and so $[e] = 1$. We know from Lemma 3.3 that the subgroup $H$ is generated by the set $\{[e] \mid e \in E^+(\Gamma)\}$. Therefore $H$ is in fact generated by the set $Y_T = \{[e] \mid e \in E^+ T\}$. It remains to show that $Y_T$ is in fact a free basis for $H$.

To see this it suffices to show that any nontrivial freely reduced word in $Y_T^{\pm 1}$ defines a nontrivial element of $F(X)$. Suppose $h = [e_1] \cdot \cdots \cdot [e_k]$ where $k \geq 1$, $e_i \in Y^+ \cup (Y^+)^{-1}$ and $e_i \in E(\hat{\Gamma} - T)$ and $e_i \neq e_{i+1}^{-1}$. We need to show that $h \neq 1$. Since $T$ is a tree, for any vertices $u, u'$ of $\Gamma$ the path-reduced form of the path $[u, v]_T [v, u']_T$ is the path $[u, u']_T$. By definition of $h$ and of $[e_i]$ we have $h = \overline{p}$ where $p$ is the following path from $v$ to $v$ in $\Gamma$:

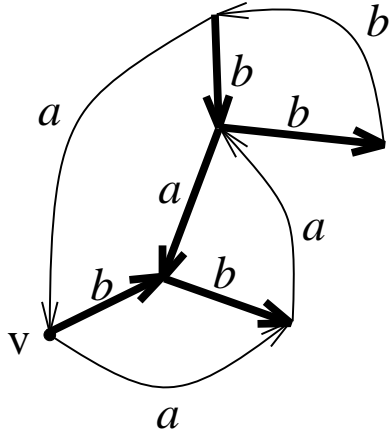

FIGURE 5. Free basis corresponding to a spanning tree. In this figure the spanning tree $T$ is highlighted. The free basis $Y_T$ for the subgroup $H = L(\Gamma, v)$ is $\{a^{-1}bab^{-1}, ba^{-1}b^3ab^{-1}, b^2a^2b^{-1}, ab^{-2}\}$.

$$p = [v, o(e_1)]_T e_1 [t(e_1), v]_T [v, o(e_2)]_T e_2 \dots [v, o(e_k)]_T e_k [t(e_k), v].$$

Also by the remark above the path $p$ can be transformed by path-reductions to the path

$$p' = [v, o(e_1)]_T e_1 [t(e_1), o(e_2)]_T e_2 \dots e_i [t(e_i), o(e_{i+1})]_T e_{i+1} \dots e_k [t(e_k), v]_T.$$

We claim that the path $p'$ is already path-reduced. Note that the segments $[t(e_i), o(e_{i+1})]_T$, $[v, o(e_1)]_T$, $[t(e_k), v]_T$ are contained in the tree $T$ and the edges $e_i$ are outside of $T$. Hence if $p'$ is not path-reduced, then for some $i$ $t(e_i) = o(e_{i+1})$ and $e_i = e_{i+1}^{-1}$. However this is impossible by our assumption on $h$.

Thus $p'$ is indeed a nontrivial reduced path in $\Gamma$. Since $\Gamma$ is folded, this implies that the label of $p'$ is a nontrivial freely reduced word in $X$ which therefore represents a nontrivial element of $F(X)$. But $p'$ was obtained by a series of path-reductions from $p$ and hence

$$1 \neq \mu(p') = \overline{\mu(p')} = \overline{\mu(p)} = h$$

and $h \neq 1$ as required. $\qquad\qquad\square$

**Proposition 6.2.** *Let $\Gamma$ be a folded connected $X$-digraph and let $\Gamma'$ be a connected subgraph of $\Gamma$. Let $v$ be a vertex of $\Gamma'$.*

*Then $H = L(\Gamma', v)$ is a free factor of $G = L(\Gamma, v)$. Further, if $\Gamma'$ does not contain $Core(\Gamma, v)$, then $H \neq G$. In particular, if both $\Gamma$ and $\Gamma'$ are core graphs with respect to $v$ and $\Gamma' \subsetneq \Gamma$ then $H \neq G$.*

*Proof.* Let $T'$ be a spanning tree in $\Gamma'$. Then there exists a spanning tree $T$ of $\Gamma$ such that $T'$ is a subgraph of $T$. By Lemma 6.1 the set $Y_T = \{[e] \,|\, e \in E^+(\Gamma - T)\}$ is a free basis of $G$. Moreover, $Y_T$ is a disjoint union of the sets:

$$Y_T = \{[e] \,|\, e \in E^+(\Gamma - T), e \in \Gamma'\} \cup \{[e] \,|\, e \in E^+(\Gamma - T), e \notin E\Gamma'\} =$$
$$= \{[e] \,|\, e \in E^+(\Gamma' - T')\} \cup \{[e] \,|\, e \in E^+(\Gamma - T), e \notin E\Gamma'\} = Y_{T'} \cup Z$$

Thus $Y_{T'}$ is a subset of $Y_T$. Since $Y_T$ is a free basis of $G$ and $Y_{T'}$ is a free basis of $G$, this implies that

$$G = F(Y_T) = F(Y_{T'}) * F(Z) = H * F(Z)$$

and $H$ is a free factor of $G$ as required.

Suppose now that $\Gamma'$ does not contain $Core(\Gamma, v)$. We claim that there is a positive edge $e$ of $\Gamma$ which does not belong to $\Gamma'$ and is not in $T$. Assume this is not the case. Then all edges outside of $\Gamma'$ lie in $T$. Hence $\Gamma - \Gamma' \subset T$ is a union of disjoint trees. This implies that $Core(\Gamma, v)$ is contained in $\Gamma'$, contrary to our assumptions. Thus the claim holds and hence $Z \neq \emptyset$. Therefore $H \neq G$ as required. $\qquad \square$

**Corollary 6.3.** *Let $\Gamma$ be a connected folded $X$-digraph which is a core graph with respect to a vertex $v$. Let $\Gamma'$ be a connected subgraph of $\Gamma$, containing $v$. Suppose that $\Gamma'$ is also a core graph with respect to $v$. Then $L(\Gamma', v)$ is a free factor of $L(\Gamma, v)$ and $L(\Gamma', v) \neq L(\Gamma, v)$.*

We will show next that the procedure described in Lemma 6.1 can be improved to obtain a Nielsen-reduced free basis of a finitely generated subgroup of $F(X)$. Here we need to recall the following important definition.

**Definition 6.4** (Nielsen set). Let $S$ be a set of nontrivial elements of the free group $F(X)$ such that $S \cap S^{-1} = \emptyset$. We say that $S$ is *Nielsen reduced with respect to the free basis $X$* if the following conditions hold:
(1) If $u, v \in S \cup S^{-1}$ and $u \neq v^{-1}$ then $|u \cdot v|_X \geq |u|_X$ and $|u \cdot v|_X \geq |v|_X$.
(2) If $u, v, w \in S \cup S^{-1}$ and $u \neq w^{-1}, v \neq w^{-1}$ then $|u \cdot w \cdot v|_X > |u|_X + |v|_X - |w|_X$.

Condition (1) means that no more than a half of $u$ and no more than a half of $v$ freely cancels in the product $u \cdot v$. Condition (2) means that at least one letter of $w$ survives after all free cancellations in the product $u \cdot w \cdot v$.

Recall that the vertex set of any connected graph $\Gamma$ is canonically endowed with an integer-valued metric $d_\Gamma$. Namely, the distance between any two vertices $u$, $u'$ is defined as the smallest length of an edge-path from $u$ to $v$ in $\Gamma$. An edge-path whose length is equal to the distance between its endpoints is said to be *geodesic in $\Gamma$*.

**Definition 6.5** (Geodesic tree). Let $\Gamma$ be a connected graph with a base-vertex $v$. A subtree tree $T$ in $\Gamma$ is said to be *geodesic relative to $v$* if $v \in T$ and for any vertex $u$ of $T$ the path $[v, u]_T$ is geodesic in $\Gamma$, that is a path of the smallest possible length in $\Gamma$ from $v$ to $u$.

It is easy to see that geodesic spanning trees always exist:

**Lemma 6.6.** *Let $\Gamma$ be a graph (whether finite or infinite) with a base-vertex $v$. Then there exists a geodesic relative to $v$ spanning tree $T$ for $\Gamma$.*

*Proof.* We will construct the tree $T$ inductively. If $\Gamma$ is finite or locally finite and recursive, our procedure will produce an actual algorithm for building $T$.

We will construct by induction a (possibly finite) sequence of nested trees $T_0 \subseteq T_1 \subseteq \dots$ in $\Gamma$ such that for each $n \geq 0$
(1) the vertex set of $T_n$ is precisely the ball $B(v, n)$ of radius $n$ centered at $v$ in the metric space $(V\Gamma, d_\Gamma)$ and
(2) the subtree $T_n$ is geodesic relative to $v$ in $\Gamma$.
**Step 0.** Put $T_0 = v$. Obviously $T_0$ is a tree and $T_0 = B(v, 0)$.
**Step $n$.** Suppose that the trees $T_0 \subseteq T_1 \subseteq \dots \subseteq T_{n-1}$ have already been constructed. Suppose further that for each $i \leq n-1$ we have $VT_i = B(v, i)$ and that $T_i$ is geodesic relative to $v$.

For each vertex $u$ of $\Gamma$ at the distance $n$ from $v$ choose an edge $e_u$ with terminus $u$ such that $d(o(e_u), v) = n - 1$. (For instance we can take $e_u$ to be the last edge of a geodesic path from $v$ to $u$ in $\Gamma$.) Put

$$T_n = T_{n-1} \cup \bigcup \{e_u \mid u \in V\Gamma, d(v, u) = n\}$$

It is easy to see that $T_n$ is a geodesic tree containing $T_{n-1}$ and that $VT_n = B(v, n)$.

Put $T = \bigcup_{n \geq 0} T_n$. Then $T$ is a geodesic spanning tree (relative to $v$) for $\Gamma$, as required. $\qquad \blacksquare$

**Proposition 6.7** (Nielsen basis). *Let $\Gamma$ be a folded $X$-digraph which is a core graph with respect to a vertex $v$ of $\Gamma$. Let $H = L(\Gamma, v) \leq F(X)$ and let $T$ be a spanning tree in $\Gamma$ which is geodesic with respect to $v$.*

*Then the set $Y_T$ is a Nielsen-reduced free basis of the subgroup $H$.*

*Proof.* Recall that $Y_T = \{[e] = \mu([v, o(e)]_T e[t(e), v]_T) \mid e \in E^+(\Gamma - T)\}$. Moreover, $[e]^{-1} = [e^{-1}] = \mu([v, o(e^{-1})]_T e[t(e^{-1}), v]_T)$ where $e \in E^+(\Gamma - T)\}$.

Note that since the tree $T$ is geodesic we have $||[v, o(e)]_T| - |[t(e), v]_T|| \leq 1$ for each $e \in E(\Gamma - T)$. Note also that the path $p_e = [v, o(e)]_T e[t(e), v]_T$ (where $e \in E(\Gamma - T)$) is path-reduced and therefore its label is freely reduced.

This means that for the path $p_e$ (with $e \in E(\Gamma - T)$) we have

$$|[v, o(e)]_T| \leq (1/2)|p_e| = (1/2)|[e]|_X \quad \text{and} \quad |[t(e), v]_T| \leq (1/2)|p_e| = (1/2)|[e]|_X$$

and

$$|e[t(e), v]_T| \geq (1/2)|p_e| = (1/2)|[e]|_X, \quad |[v, o(e)]_T e| \geq 1/2)|p_e| = (1/2)|[e]|_X.$$

We already know from Lemma 6.1 that the set $Y_T$ is a free basis of $H$. It remains to check that $Y_T$ is Nielsen-reduced.

We first check condition (1) of Definition 6.4. Suppose $e, f \in E(\Gamma - T)$ and $e \neq f^{-1}$. Then the path-reduced form of the path $p_e p_f = [v, o(e)]_T e[t(e), v]_T[v, o(f)]_T f[t(f), v]_T$ is the path

$$p' = [v, o(e)]_T e[t(e), o(f)]_T f[t(f), v]_T.$$

The label of this path is freely reduced (since $\Gamma$ is folded and therefore it is equal to $[e] \cdot [f] \in H$. Since $p'$ has a subpath $[v, o(e)]_T e$, we conclude that $|[e] \cdot [f]|_X \geq (1/2)|[e]|_X$. Since $p'$ has a subpath $f[t(f), v]_T$, we conclude that $|[e] \cdot [f]|_X \geq (1/2)|[f]|_X$. Thus condition (1) holds.

We will now verify condition (2) of Definition. Suppose $e, f, g \in E(\Gamma - T)$ are such that $e \neq f^{-1}, g \neq f^{-1}$. We need to show that $|[e] \cdot [f] \cdot [g]|_X > |[e]|_X + |[g]|_X - |[f]|_X$.

The path-reduced form of the concatenation

$$p_e p_f p_g = [v, o(e)]_T e[t(e), v]_T[v, o(f)]_T f[t(f), v]_T[v, o(g)]_T g[t(g), v]_T$$

is the path

$$p'' = [v, o(e)]_T e[t(e), o(f)]_T f[t(f), o(g)]_T g[t(g), v]_T.$$

When $p_e p_f p_g$ is transformed to $p'$, we have to cancel the terminal segment of $p_e$ of length at most $|[v, o(f)]_T|$ and the initial segment of $p_g$ of length at most $|[t(f), v]_T|$. Since the edge $f$ of $p_f$ survives, we have

$$|p''| \geq 1 + |p_e| + |p_g| - |[v, o(f)]_T| - |[t(f), v]_T| \geq 1 + |p_e| + |p_g| - (1/2)|p_f| - (1/2)|p_f| > |p_e| + |p_g| - |p_f|$$

and so

$$|[e] \cdot [f] \cdot [g]|_X > |[e]|_X + |[g]|_X - |[f]|_X$$

as required. $\qquad\square$

**Corollary 6.8.** *Let $H \leq F(X)$ be a subgroup of $F(X)$. Then $H$ is free. Moreover, there is a Nielsen-reduced free basis for $H$.*

*Proof.* Consider the graph $\Gamma(H)$. By Lemma 6.6 there exists a geodesic spanning tree $T$ in $\Gamma(H)$. Hence by Proposition 6.7 the set $Y_T$ is a Nielsen-reduced free basis for $H$. In particular $H$ is free. $\qquad\square$

**Historical Note 6.9.** The fact that finitely generated subgroups in a free group $F(X)$ of finite rank are free must have been known to Klein and Fricke in 1880's. The first formal proof is due to Jacob Nielsen and appeared in 1923 (see [38]). Later Nielsen generalized it to arbitrary subgroups of $F(X)$ in 1955 [36]. Today this statement appears trivial by topological considerations. Indeed, $F(X)$ can be realized as the fundamental group of a topological graph. A covering space of a graph is a graph and fundamental groups of graphs are free.

However, the existence of a *Nielsen-reduced* free basis is not at all trivial, even for finitely generated subgroups of $F(X)$. The standard argument, due to J.Nielsen himself, uses a rather complicated machinery of Nielsen transformations and a particular choice of well-ordering on $F(X)$. The advantage of our approach is that we use only very elementary combinatorial objects.

**Remark 6.10.** It can be shown, for example by looking at abelianizations, that for a free group any two bases have the same cardinality. If $F$ is a finitely generated free group, then the cardinality of any free basis of $F$ is called *the rank of* $F$ and denoted $rk(F)$. It is easy to see that for finitely generated free groups $F$ and $F_1$ the free product $F * F_1$ is also a free group and $rk(F * F_1) = rk(F) + rk(F_1)$.

## 7. Basic algebraic and algorithmic properties of subgroup graphs

In this section we will restrict our attention to finitely generated subgroups of free groups. Moreover, we will be particularly interested in algorithmic aspects of various group-theoretic statements.

Till the end of this section let $X$ be a finite set and let $F = F(X)$ be the free group on $X$.

**Proposition 7.1** (Constructing $\Gamma(H)$). *There is an algorithm which, given finitely many freely reduced $X$-words $h_1, \ldots, h_k$, constructs the graph $\Gamma(H)$, where $H = \langle h_1, \ldots, h_k \rangle$.*

*Proof.* This follows immediately from the proof of Proposition 3.8. Indeed, if we start with a wedge of $k$ circles with the words $h_1, \ldots, h_k$ written on them and then perform all possible foldings, the resulting $X$-digraph will be $\Gamma(H)$. Note that the original wedge-of-circles graph has $M = |h_1|_X + \cdots + |h_k|_X$ edges. Since every folding reduces the number of edges by one, we will obtain the graph $\Gamma(H)$ in at most $M$ steps. $\qquad\square$

**Proposition 7.2** (Generalized word problem). *The group $F = F(X)$ has solvable generalized word problem.*

*That is, there exists an algorithm which, given finitely many freely reduced $X$-words $h_1, \ldots, h_k$ and given a freely reduced $X$-word $g$ decides whether or not $g$ belongs to the subgroup $H = \langle h_1, \ldots, h_n \rangle$ of $F$.*

*Proof.* We first construct the graph $\Gamma(H)$ as it is done in Proposition 7.1. Recall that $L(\Gamma(H), 1_H) = H$ by definition of $\Gamma(H)$.

Thus to decide if $g \in H$, we simply have to check whether or not the graph $\Gamma(H)$ accepts the word $g$. Recall that $\Gamma(H)$ is folded and that the language $L(\Gamma(H), 1_H)$ consists of all labels of reduced paths in $\Gamma(H)$ from $1_H$ to $1_H$. Given a freely reduced word $g = x_1 \ldots x_m$, $x_i \in \Sigma = X \cup X^{-1}$ we can check whether $g \in L(\Gamma(H), 1_H)$ as follows.

First we check if there is an edge in $\hat{\Gamma}(H)$ with origin $1_H$ and label $x_1$ (since $\Gamma$ is folded, such an edge is unique if it exists). If yes, we move to the terminal vertex of this edge, which we denote $v_1$. If no, we terminate the process. We then check if there is an edge in $\hat{\Gamma}(H)$ with origin $v_1$ and label $x_2$. If yes, we move to the terminal vertex of this edge, which we denote $v_2$. If no, we terminate the process. By repeating this procedure at most $k$ times we either will terminate the process, in which case we conclude that $g \notin L(\Gamma(H), 1_H)$, or we find a vertex $v_k$ such that the word $g$ is the label of a reduced path in $\Gamma(H)$ from $1_H$ to $v_k$. If $v_k \neq 1_H$, we conclude that $g \notin L(\Gamma(H), 1_H)$ and if $v_k = 1_H$, we conclude that $g \in L(\Gamma(H), 1_H)$. $\qquad\square$

Recall that a core graph is connected and has at most one vertex of degree one (namely the base-vertex).

**Definition 7.3** (Type of a core graph). Let $\Gamma$ be a folded $X$-digraph which is a core graph with respect to some vertex. Suppose that $\Gamma$ has at least one edge.

If every vertex of $\Gamma$ has degree at least two, we set the *type* of $\Gamma$, denoted $Type(\Gamma)$ to be equal to $\Gamma$.

Suppose now that $\Gamma$ has a vertex $v$ of degree one (such a vertex is unique since $\Gamma$ is a core graph). Then there exists a unique vertex $v'$ of $\Gamma$ with the following properties:

(a) There is a unique $\Gamma$-geodesic path $[v, v']$ from $v$ to $v'$ and every vertex of this geodesic, other than $v$ and $v'$, has degree two.

(b) The vertex $v'$ has degree at least three.

Let $\Gamma'$ be the graph obtained by removing from $\Gamma$ all the edges of $[v, v']$ and all the vertices of $[v, v']$ except for $v'$. Then $\Gamma'$ is called the *type* of $\Gamma$ and denoted $Type(\Gamma)$.

Finally, if $\Gamma$ consists of a single vertex, we set $Type(\Gamma) = \Gamma$.

**Remark 7.4.** Observe that $Type(\Gamma)$ is an $X$-digraph which does not have a distinguished base-vertex. Note also that the graph $\Gamma' = Type(\Gamma)$ is a core-graph with respect to any of its vertices. Moreover in the graph $\Gamma'$ every vertex has degree at least two and $\Gamma' = Type(\Gamma')$.

**Lemma 7.5.** *Let $\Gamma$ be a folded core graph (with respect to one of its vertices). Let $v$ and $u$ be two vertices of $\Gamma$ and let $q$ be a reduced path in $\Gamma$ from $v$ to $u$ with label $g \in F(X)$. Let $H = L(\Gamma, v)$ and $K = L(\Gamma, u)$.*

*Then $H = gKg^{-1}$.*

*Proof.* Let $p$ be a reduced path in $\Gamma$ from $u$ to $u$ (so that $\mu(p) = k \in K$). Then the path $p' = qpq^{-1}$ is a path from $v$ to $v$ with label $\mu(p') = \mu(q)\mu(p)(\mu(q))^{-1}$. Thus the freely reduced form of $\mu(p')$ is equal to the element $g \cdot k \cdot g^{-1} \in F(X)$. The path $p'$ can be transformed by several path-reductions to a reduced path $p''$ from $v$ to $v$. Hence $\mu(p'') \in L(\Gamma, v) = H$. On the other hand $\mu(p'') = \overline{\mu(p')} = g \cdot k \cdot g^{-1}$. Thus we have shown that for any $k \in K$ we have $g \cdot k \cdot g^{-1} \in H$ so that $gKg^{-1} \subseteq H$. A symmetric argument shows that $g^{-1}Hg \subseteq K$ and thus $H \subseteq gKg^{-1}$. Therefore $H = gKg^{-1}$, as required. $\square$

**Lemma 7.6.** *Let $H \le F(X)$ and let $\Gamma = \Gamma(H)$.*

*Let $g \in F(X)$ be a nontrivial freely reduced word in $X$. Let $g = yz$ where $z$ is the maximal terminal segment of the word $g$ such that there is a path with label $z^{-1}$ in $\Gamma$ starting at $1_H$ (such a path is unique since $\Gamma$ is folded). Denote the end-vertex of this path by $u$. Let $\Delta'$ be the graph obtained from $\Gamma$ as follows. We attach to $\Gamma$ at $u$ the segment consisting of $|y|$ edges with label $y^{-1}$, as read from $u$. Let $u'$ be the other end of this segment. Put $\Delta'' = Core(\Gamma, u')$.*

*Then $(\Delta'', u') = (\Gamma(K), 1_K)$ where $K = gHg^{-1}$.*

*Proof.* Let $\Delta$ be the graph obtained by attaching to $\Gamma$ at $1_H$ a segment of $|g|_X$ edges, labeled $g$, which has origin $v$ and terminus $1_H$. It is easy to see that $\overline{(L(\Delta, v))} = gHg^{-1}$. However, the graph $\Delta$ is not necessarily folded. If we fold the terminal segment $z$ of $g$ (as defined in the statement of Lemma 7.6) onto the path in $\Gamma$ with origin $1_H$ and label $z^{-1}$, the resulting graph $\Delta'$ is obviously folded. Therefore $L(\Delta', u') = H$ (where $u'$ is the image of $v$ in $\Delta'$). It can still happen that $\Delta'$ is not a core graph with respect to $u'$. (see Figure). However the graph $\Delta'' = Core(\Delta', u')$ is a folded, is a core graph with respect to $u'$ and has the property $L(\Delta'', u') = K$. Hence $(\Delta'', u') = (\Gamma(K), 1_K)$, as required. $\square$

The following observation will allow us to decide when two finitely generated subgroups of $F(X)$ are conjugate.

**Proposition 7.7.** *[Conjugate subgroups] Let $H$ and $K$ be subgroups of $F(X)$.*

*Then $H$ is conjugate to $K$ in $F(X)$ if and only if the graphs $Type(\Gamma(H))$ and $Type(\Gamma(K))$ are isomorphic as $X$-digraphs.*

*Proof.* Suppose that $Type(\Gamma(H)) = Type(\Gamma(K)) = \Gamma$. Let $v$ be a vertex of $\Gamma$. It follows from Lemma 7.5 that the subgroup $L(\Gamma, v)$ is conjugate to both $H$ and $K$, so that $H$ is conjugate to $K$.

Suppose now that $K$ is conjugate to $H$, that is $K = gHg^{-1}$ for some $g \in F(X)$. Lemma 7.6 implies that $Type(\Gamma(H)) = Type(\Gamma(K))$, as required. $\square$

**Corollary 7.8** (Conjugacy problem for subgroups of free groups)**.** *There is an algorithm which, given finitely many freely reduced words $h_1, \ldots, h_s, k_1, \ldots, k_m$, decides whether or not the subgroups $H = \langle h_1, \ldots, h_s \rangle$ and $K = \langle k_1, \ldots, k_m \rangle$ are conjugate in $F(X)$*

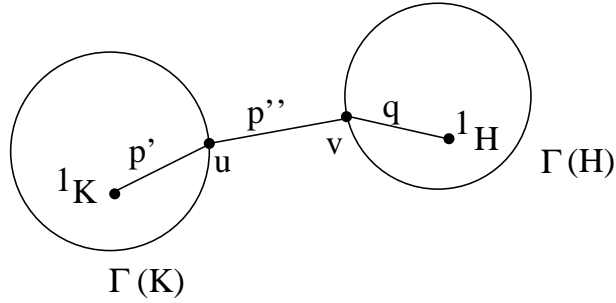

FIGURE 6

*Proof.* First we construct the graphs $\Gamma(H)$ and $\Gamma(K)$, as described in Proposition 7.1. By chopping-off segments with origin $1_H$, $1_K$ if necessary, we then obtain the graphs $\Gamma = Type(\Gamma(H))$ and $\Delta = Type(\Gamma(K))$. It remains to check whether $\Gamma$ and $\Delta$ are isomorphic as $X$-digraphs. This can be done by fixing a vertex $v \in V\Gamma$ and comparing for each vertex $u \in V\Delta$ the based $X$-digraphs $(\Gamma, v)$ and $(\Delta, u)$. □

**Proposition 7.9.** *Let $H, K \leq F(X)$ and $g \in F(X)$. Let $g = yz$ where $z$ is the maximal terminal segment of $g$ such that $z^{-1}$ is the label of a path in $\Gamma(H)$ with origin $1_H$. Let $y = y'y''$ where $y'$ is the maximal initial segment of $y$ which is the label of a path in $\Gamma(K)$ with origin $1_K$.*

*Suppose that the word $y''$ is nontrivial. Then*

$$\langle gHg^{-1}, K \rangle = gHg^{-1} * K.$$

*Proof.* Consider the graph $\Gamma$ obtained in the following way. First we take an edge-path $p_g$ of length $|g|_X$ with label $g$ and attach its origin to $1_K$ and its terminus to $1_H$. The path $p_g$ is a concatenation of the form $p_g = p'p''q$ where $\mu(p') = y'$, $\mu(p'') = y''$ and $\mu(q) = z$.

Then we fold the segment $q$ of $p_g$ onto the path (starting at $1_H$) with label $z^{-1}$ in $\Gamma(H)$ and we fold the segment $y'$ of $p_g$ onto the path (starting at $1_K$) in $\Gamma(K)$. We denote the resulting graph $\Gamma$.

Thus $\Gamma$ consists of the graphs $\Gamma(H)$ and $\Gamma(K)$ joined by a segment with label $y''$ (see Figure 6). By a slight abuse of notation we will still call this segment $p''$ and we will also refer to the images of $\Gamma(H)$ and $\Gamma(K)$ in $\Gamma$ as $\Gamma(H)$ and $\Gamma(K)$. The construction of $\Gamma$ implies that it is a folded $X$-digraph.

It is not hard to see now that $L(\Gamma, 1_K) = gHg^{-1} * K$. Indeed, choose a spanning tree $T_H$ is $\Gamma(H)$ and a spanning tree $T_K$ in $\Gamma(K)$. Then $T = T_H \cup T_K \cup p''$ is a spanning tree in $\Gamma$. We can use it to produce a free basis $Y_T$ for $L(\Gamma, 1_K)$ as described in Lemma 6.1. Then $Y_T$ is a disjoint union of two subsets $Y_T = Y_H \sqcup Y_K$ where $Y_H$ are the generators corresponding to the edges in $\Gamma(H) - T_H$ and $Y_K$ are the generators corresponding to the edges in $\Gamma(K) - T_K$.

Each $f \in Y_H$ is the label of a path $\alpha_f = [1_K, u]_{T_K} p''[v, o(e)]_{T(H)} e[t(e), v]_{T_H} (p'')^{-1}[u, 1_K]_{T_K}$ where $e \in E^+(\Gamma(H) - T_H)$. The path $p_f$ can be obtained by path-reductions from

$$\alpha'_f = ([1_K, u]_{T_K}(p')^{-1})p'p''q(q^{-1}[v, 1_H]_{T_H})[1_H, o(e)]_{T_H} e[t(e), 1_H]_T ([1_H, v]_{T_H} q)q^{-1}(p'')^{-1}(p')^{-1}(p'[u, 1_K]_{T_K})$$

that is

$$\alpha'_f = ([1_K, u]_{T_K}(p')^{-1})p_g(q^{-1}[v, 1_H]_{T_H})p_e([1_H, v]_{T_H} q)p_g^{-1}(p'[u, 1_K]_{T_K})$$

where $p_e = [1_H, o(e)]_{T_H} e[t(e), 1_H]_T$ is the path from $1_H$ to $1_H$ in $\Gamma(H)$ defining one of the standard free generators of $H$ relative the tree $T_H$.

Since $[1_K, u]_{T_K}(p')^{-1}$ is a path in $\Gamma(K)$ from $1_K$ to $1_K$, the label of this path determines an element $k \in K$. Similarly, since $q^{-1}[v, 1_H]_{T_H} = ([1_H, v]_{T_H} q)^{-1}$ is a path in $\Gamma(H)$ from $1_H$ to $1_H$, the freely reduced form of its label is an element $h \in H$.

Since the set $Y_{T_H} = \{\mu(p_e) \mid e \in E^+(\Gamma(H) - T_H)\}$ is a free basis of $H$, the set $Z_H = g \cdot h \cdot Y_{T_H} \cdot h^{-1} \cdot g^{-1}$ is a free basis of $gHg^{-1}$. Thus $Y_H = \overline{k \cdot Z_H k^{-1}}$.

On the other hand the set $Y_K$ is clearly exactly the free basis of $K$ corresponding to the spanning tree $T_K$ in $\Gamma(K)$.

Since we know that $Y_T = Y_H \cup Y_K$ is a free basis of the subgroup $L(\Gamma, 1_K)$ and $F(Y_K) = K$ it follows that $(k^{-1}Y_H k) \cup Y_T = Z_H \cup Y_T$ is also a free basis of $L(\Gamma, 1_K)$. Therefore

$$L(\Gamma, 1_K) = F(Z_H \cup Y_K) = F(Z_H) * F(Y_K) = gHg^{-1} * K$$

as required. $\qquad\square$

**Proposition 7.10.** *Let $H \leq F(X)$ be a finitely generated subgroup and let $\Gamma = \Gamma(H)$. Let $g \in F(X)$ be such that $g^n \in H$ for some $n \geq 1$. Then there exists $m \leq \#V\Gamma$, $m > 0$, such that $g^m \in H$.*

*Proof.* We may assume that $g \notin H$, since otherwise the statement is obvious.

Let $g = f d f^{-1}$ where the subword $f$ of $g$ is cyclically reduced. Then for any $i \geq 1$ the freely reduced form of $g^i$ is $f d^i f^{-1}$.

Let $n \geq 1$ be the smallest positive integer such that $g^n \in H$. Since $g \notin H$, we have $n > 1$. Suppose that $n > \#V\Gamma$.

Since $g^n \in H$, the freely reduced word $f d^n f^{-1}$ is the label of a reduced path $p$ in $\Gamma$ from $1_H$ to $1_H$. For each $k = 0, 1, \ldots, n$ let $v_k$ be the end-vertex of the initial segment of $p$ with label $f d^k$. Since $n > \#V\Gamma$, for some $1 \leq i < j$ we have $v_i = v_j = v$. Then the subpath $\alpha$ of $p$ from $v_i$ to $v_j$ has label $d^{j-i}$, so that $p = p_1 \alpha p_2$. Put $p'$ be the path obtained from $p$ by removing the subpath $\alpha$, that is $p' = p_1 p_2$. It is obvious that the label of $p'$ is $f d^{n-(j-i)} f^{-1}$. Moreover, since $i \geq 1$, we have $1 \leq n - (j-i) < n$, $n - (j-i) > 0$ and so the word $f d^{n-(j-i)} f^{-1}$ is freely reduced. Since $\Gamma$ is folded, the path $p'$ is reduced as well.

Therefore by definition of $\Gamma(H)$ we have $g^{n-(j-i)} \in H$ and $1 \leq n - (j-i) < n$. This contradicts the choice of $n$. $\qquad\square$

**Corollary 7.11.** *There exists an algorithm which, given finitely many freely reduced words*

$$g, h_1, \ldots, h_s \in F(X)$$

*decides whether or not some nonzero power of $g$ belongs to the subgroup $H = \langle h_1, \ldots, h_s \rangle$.*

The use of subgroup graphs also allows us to easily determine whether one subgroup is conjugate to a subgroup of another.

**Lemma 7.12.** *Let $K, H$ be subgroups of $F(X)$. Then there is $g \in F(X)$ with $gKg^{-1} \leq H$ if and only if there exists a morphism of (non-based) $X$-digraphs $\pi : Type(\Gamma(K)) \longrightarrow Type(\Gamma(H))$.*

*Proof.* Suppose $K$ is conjugate to a subgroup of $H$. Since conjugation preserves type, we can assume that $\Gamma(H) = Type(\Gamma(H))$. Let $g \in F(X)$ be such that $gKg^{-1} \leq H$. By Proposition 4.3 there exists a morphism $\beta : \Gamma(gKg^{-1}) \longrightarrow \Gamma(H)$. However by Lemma 7.6 $Type(\Gamma(K)) = Type(\Gamma(gKg^{-1}))$ is a subgraph of $\Gamma(gKg^{-1})$. The restriction of $\beta$ to this subgraph produces the desired morphism $\pi : Type(\Gamma(K)) \longrightarrow Type(\Gamma(H))$.

Suppose now that a morphism $\pi : Type(\Gamma(K)) \longrightarrow Type(\Gamma(H))$ exists. Let $f$ be the label of the shortest path in $\Gamma(K)$ from $1_K$ to $Type(\Gamma(K))$ and let $u$ be the terminal vertex of this path. Let $v = \pi(u)$ be the image of this vertex in $\Gamma(H)$.

Let $c$ be the label of a shortest path from $1_H$ to $v$ in $\Gamma(H)$. Thus by Lemma 7.5 $L(Type(K), u) = f^{-1}Kf$ and $L(Type(H), v) = c^{-1}Hc$, so that $(Type(K), u) = (\Gamma(f^{-1}Kf), 1_{f^{-1}Kf})$ and $(Type(H), v) = (\Gamma(c^{-1}Hc), 1_{c^{-1}Hc})$. Since $\pi$ can also be considered as a morphism $\pi : (\Gamma(f^{-1}Kf), 1_{f^{-1}Kf}) \longrightarrow (\Gamma(c^{-1}Hc), 1_{c^{-1}Hc})$, Proposition 4.3 implies that

$$f^{-1}Kf \leq c^{-1}Hc, \text{ and } K \leq (f \cdot c^{-1})H(f \cdot c^{-1})^{-1},$$

as required. Then $K$ is conjugate to a subgroup of $H$ by Proposition 7.7 and Proposition 4.3. $\qquad\square$

**Corollary 7.13.** *There exists an algorithm which, given finitely many freely reduced words*

$$k_1, \ldots, k_s, h_1, \ldots, h_m \in F(X),$$

*decides whether or not there is $g \in F(X)$ such that $gKg^{-1} \leq H$, where $K = \langle k_1, \ldots, k_s \rangle$, $H = \langle h_1, \ldots, h_m \rangle$. Moreover, the algorithm will produce one such $g$ if it exists.*

*Proof.* We first construct the graphs $\Gamma(K), \Gamma(H)$ and their type-subgraphs $Type(\Gamma(K))$, $Type(\Gamma(H))$. Let $u$ be the vertex of $Type(\Gamma(K))$ closest to $1_K$ (thus $u = 1_K$ if $\Gamma(K) = Type(\Gamma(K))$ and $u$ is the terminal vertex of the segment from $1_K$ to $Type(\Gamma(K))$ otherwise. Let $f$ be the label of the unique reduced path from $1_K$ to $u$ in $\Gamma(K)$.

We then enumerate the vertices of $Type(\Gamma(H))$ and for each vertex $v$ check whether there exists a morphism $\pi : Type(\Gamma(K)) \longrightarrow Type(\Gamma(H))$ which takes $u$ to $v$. If no such vertex of $Type(\Gamma(H))$ is found, $K$ is not conjugate to a subgroup of $H$ by Lemma 7.12. If such $v$ is found, then by the proof of Lemma 7.12 $gKg^{-1} \leq H$ where $g = f \cdot c^{-1}$ and $c$ is the label of a reduced path from $1_H$ to $v$ in $\Gamma(H)$. $\square$

## 8. Finite Index, Commensurability and Marshall Hall Theorem

**Definition 8.1** (Regular digraphs). An $X$-digraph $\Gamma$ is said to be *$X$-regular* if for every vertex $v$ of $\Gamma$ and every $x \in X \cup X^{-1}$ there is exactly one edge in $\hat{\Gamma}$ with origin $v$ and label $x$.

We need to prove the following simple graph-theoretic statement.

**Lemma 8.2.** *Let $\Gamma$ be a finite connected digraph. Let $T$ be a spanning tree in $\Gamma$. Then $\#E^+(\Gamma - T) = \#E^+\Gamma - \#V\Gamma + 1$.*

*In particular, if $\Gamma$ is a finite connected folded $X$-digraph then $\#Y_T = \#E^+\Gamma - \#V\Gamma + 1$. Thus if $v$ is a vertex of $\Gamma$ and $H = L(\Gamma, v)$ then $rk(H) = \#E^+\Gamma - \#V\Gamma + 1$.*

*Proof.* It is easy to prove by induction on the number of edges that for a finite tree the number of vertices minus the number of positive edges is equal to one.

Hence $\#E^+T = \#VT - 1 = \#V\Gamma - 1$. Therefore

$$\#E^+\Gamma = \#E^+T + \#E^+(\Gamma - T) = \#V\Gamma - 1 + \#E^+(\Gamma - T)$$

and so

$$\#E^+(\Gamma - T) = \#E^+\Gamma - \#V\Gamma + 1$$

as required. $\square$

The above statement can be used to compute the rank of a finitely generated subgroup $H \leq F(X)$ via its graph $\Gamma(H)$.

We can now give a criterion for a subgroup $H \leq F(X)$ to be of finite index in $F(X)$.

**Proposition 8.3** (Finite index subgroups). *Let $H$ be a subgroup of $F(X)$.*

*Then $|F(X) : H| < \infty$ if and only if $\Gamma(H)$ is a finite $X$-regular graph. In this case $|F(X) : H| = \#V\Gamma(H)$.*

*Proof.* Assume first that $|F(X) : H| = s < \infty$. Then $H$ is finitely generated and thus $\Gamma(H)$ is finite. Suppose that $\Gamma(H)$ is not $X$-regular. Then there is a vertex $v$ of $\Gamma$ and a letter $x \in X \cup X^{-1}$ such that there is no edge labeled $x$ with origin $v$ in $\hat{\Gamma}$. Choose a reduced path $p$ from $1_H$ to $v$ in $\Gamma$ and let $w = \mu(p)$. Then the word $wx$ is freely reduced and, moreover, no word with initial segment $wx$ is accepted by $\Gamma(H)$. Thus for any freely reduced word $y$ whose first letter is not $x^{-1}$, we have $wxy \notin H$.

Since $|F(X) : H| = s < \infty$, there exist $s$ elements $g_1, \ldots, g_s \in F(X)$ such that

$$F(X) = Hg_1 \cup \cdots \cup Hg_s$$

Let $M = \max\{|g_i|_X \mid i = 1, \ldots, s\}$ and let $f = wx^{M+1}$. Then for some $g_i$ we have $f \cdot g_i^{-1} \in H$. By the choice of $M$ the freely reduced form $z$ of $f \cdot g_i^{-1}$ has initial segment $wx$ and therefore $z \notin L(\Gamma(H), 1_H) = H$. This contradicts our assumption that $f \cdot g_i^{-1} \in H$.

Suppose now that $\Gamma(H)$ is an $X$-regular finite graph. We want to show that the index of $H$ in $F(X)$ is finite and, moreover, that it is equal to the number of vertices in $\Gamma(H)$.

For each vertex $v$ of $\Gamma(H)$ choose a reduced path $p_v$ from $1_H$ to $v$ in $\Gamma(H)$ and denote $g_v = \mu(p_v)$. We claim that

$$(1) \qquad\qquad F(X) = \cup\{H \cdot g_v \mid v \in V\Gamma(H)\}.$$

Let $f \in F(X)$ be an arbitrary freely reduced word. Since $\Gamma(H)$ is $X$-complete, there is a path $\alpha$ in $\Gamma(H)$ with origin $1_H$ and label $f$. Let $u$ be the terminal vertex of this path. Then $\alpha p_u^{-1}$ is a path in $\Gamma(H)$ from $1_H$ to $1_H$ with label $fg_u^{-1}$. Therefore the freely reduced form of the word $fg_u^{-1}$ is accepted by $\Gamma(H)$ and so $f \cdot g_u^{-1} \in H, f \in H \cdot g_u$. Since $f \in F(X)$ was chosen arbitrarily, we have verified that (1) holds.

This already implies that $H$ has finite index in $F(X)$. To prove the proposition it remains to check that if $v$ and $u$ are distinct vertices of $\Gamma(H)$ then $Hg_v \neq Hg_u$. Suppose this is not so and for some $v \neq u$ we have $Hg_v = Hg_u$. Then $g_v \cdot g_u^{-1} = h \in H$. The path $p = p_v p_u^{-1}$ has label $g_v g_u^{-1}$ and hence the path-reduced form $p'$ of $p$ has label $h$. Since $p'$ is a reduced path starting at $1_H$ with label $h \in H$, the terminal vertex of $p'$ is also $1_H$. Path-reductions do not change the end-points of a path. Therefore the terminal vertex of $p$ is $1_H$, as well as the initial vertex of $p$. However this implies that the paths $p_v$ and $p_u$ have the same terminal vertices and $v = u$, contrary to our assumptions. □

**Corollary 8.4** (Schrier's Formula). *Let $H$ be a subgroup of finite index $i$ in a finitely generated free group $F = F(X)$. Then*
$$rk(H) - 1 = i(rk(F) - 1)$$

*Proof.* Let $\Gamma = \Gamma(H)$. Then by Lemma 8.2 $rk(H) - 1 = \#E^+\Gamma - \#V\Gamma$.

Since $|F : H| = i < \infty$, the graph $\Gamma$ is $X$-regular and the degree of each vertex is $2\#X = 2rk(F)$. If we add up the degrees of all vertices of $\Gamma$, we will count every edge twice. Therefore $2rk(F)\#V\Gamma = 2\#E^+\Gamma$ and $rk(F)\#V\Gamma = \#E^+\Gamma$. Recall also that by Proposition 8.3 $\#V\Gamma = i$. Hence
$$rk(H) - 1 = \#E^+\Gamma - \#V\Gamma = rk(F)\#V\Gamma - \#V\Gamma = rk(F)i - i = i(rk(F) - 1)$$
as required. ∎

**Corollary 8.5.** *There exists an algorithm which, given finitely many freely reduced words $h_1, \ldots, h_s \in F(X)$, computes the index of the subgroup $H = \langle h_1, \ldots, h_s \rangle$ in $F(X)$.*

*Proof.* We first constrict the finite graph $\Gamma(H)$ and then check if it is $X$-regular. If no, the subgroup $H$ has infinite index in $G$ (by Proposition 8.3). If yes, the index of $H$ in $G$ is equal to the number of vertices of $\Gamma(H)$ (also by Proposition 8.3). ∎

**Definition 8.6.** Let $H$ be a subgroup of a group $G$. The *commensurator $Comm_G(H)$* of $H$ in $G$ is defined as

$$Comm_G(H) = \{g \in G \mid |H : H \cap gHg^{-1}| < \infty \text{ and } |gHg^{-1} : H \cap gHg^{-1}| < \infty\}$$

It is easy to see that $Comm_G(H)$ is a subgroup of $G$ containing $H$.

**Lemma 8.7.** *Let $H \leq F(X)$ be a nontrivial finitely generated subgroup. Then $|F(X) : H| < \infty$ if and only if $F(X) = Comm_{F(X)}(H)$.*

*Proof.* It is obvious that $|F(X) : H| < \infty$ implies $F(X) = Comm_{F(X)}(H)$. Suppose now that $H$ is a finitely generated subgroup of $F(X)$ and $F(X) = Comm_{F(X)}(H)$. Assume that $|F(X) : H| = \infty$.

Then the graph $\Gamma(H)$ is not $X$-regular. Thus there is a vertex $v \in \Gamma$ and a letter $x \in X \cup X^{-1}$ such that there is no edge labeled $x$ with origin $v$ in $\Gamma(H)$. Since $H$ is a nontrivial subgroup of infinite index in $F(X)$, the rank of $F(X)$ is at least two and so $\#(X) \geq 2$. Let $a \in X$ be a letter such that $a \neq x^{\pm 1}$.

Since for any $g \in F(X)$ we have $|gHg^{-1} : H \cap gHg^{-1}| < \infty$, for any element of $gHg^{-1}$ some power of this element belongs to $H \cap gHg^{-1}$ and so to $H$. Hence for any $g \in F(X)$ and any $h \in H$ there is $n \geq 1$ such that $g^{-1}h^n g \in H$.

Let $h \in H$ be a nontrivial element, so that $h$ is a freely reduced word in $X$. Let $y \in X \cup X^{-1}$ be the first letter of $h$ and let $z \in X \cup X^{-1}$ be the last letter of $h$. Then for any $m \geq 1$ the freely reduced form of $h^n$ begins with $y$ and ends with $z$.

Let $w$ be the label of a reduced path in $\Gamma(H)$ from $1_H$ to $v$. Since there is no edge labeled $x$ with origin $v$ in $\Gamma(H)$, any freely reduced word with initial segment $wx$ does not belong to $H$.

Put $q = y$ if $y \neq z^{-1}$. If $y = z^{-1}$ and $y \in \{x, x^{-1}\}$, put $q = a$. If $y = z^{-1}$ and $y \notin \{x, x^{-1}\}$, put $q = x$. Then for any $m \geq 1$ the word $q\overline{h^m}q^{-1}$ is freely reduced. (Recall that $\overline{h^m}$ is the freely reduced form of $h^m$.)

Choose a freely reduced word $w'$ such that the word $wxw'q$ is freely reduced. This is obviously possible since $X$ has at least two elements. Put $g = wxw'q$. By our assumptions there is $n \geq 1$ such that the element $g \cdot h^n g^{-1}$ belongs to $H$. However by the choice of $q$ the freely reduced form of $gh^n g^{-1}$ is $wxw'qy \ldots zq^{-1}(w')^{-1}x^{-1}w^{-1}$. The word $wxw'qy \ldots zq^{-1}(w')^{-1}x^{-1}w^{-1}$ has initial segment $wx$ and hence cannot represent an element of $H$. This yields a contradiction. □

**Corollary 8.8** (Greenberg-Stallings Theorem). [19], [43] *Let $H, K$ be finitely generated subgroups of $F(X)$ such that $H \cap K$ has finite index in both $H$ and $K$.*

*Then $H \cap K$ has finite index in the subgroup $\langle H \cup K \rangle$.*

*Proof.* Let $G = \langle H \cup K \rangle$ and $L = H \cap K$. Since $|H : L| < \infty$ and $|K : L| < \infty$, we have $H \leq Comm_{F(X)}(L)$, $K \leq Comm_{F(X)}(L)$ and hence $G = \langle H \cup K \rangle \leq Comm_{F(X)}(L)$.

Both $H$ and $K$ are finitely generated and therefore so are $G$ and $L$ (recall that $L$ is a subgroup of finite index in a finitely generated free group $H$). Hence by Corollary 6.8 $G$ is a free group of finite rank, $G = F(S)$ for some finite set $S$. Since $G \leq Comm_{F(X)}(L)$, we have $G = Comm_G(L)$. This, by Lemma 8.7, implies that $|G : L| < \infty$, as required. □

**Proposition 8.9.** *Let $H$ be a nontrivial finitely generated subgroup of $F(X)$. Then $H$ has finite index in $Comm_{F(X)}(H)$.*

*Proof.* Let $G = Comm_{F(X)}(H)$. If $G$ is finitely generated then $G$ is a free group of finite rank containing $H$ and $G = Comm_G(H)$. Hence $|G : H| < \infty$ by Lemma 8.7.

Suppose now that $G$ is not finitely generated. Then $G = F(Y)$ where $Y$ is an infinite free basis of $G$. Since $H$ is finitely generated, the elements of a finite generating set for $H$ involve only finitely many letters of $Y$. Thus $H$ is contained in a finitely generated free factor $A$ of $G$ and $G = A * B$ where $B$ is a free group of infinite rank. However $G = Comm_G(H)$ and so for each $g \in G$ we have $|H : gHg^{-1} \cap H| < \infty$. However for any $b \in B$ $bAb^{-1} \cap A = 1$ and so $bHb^{-1} \cap H = 1$. This contradicts the fact that $H$ is nontrivial and thus infinite. □

**Lemma 8.10.** *Let $\Gamma$ be a finite folded $X$-digraph. Then there exists a finite folded $X$-regular digraph $\Gamma'$ such that $\Gamma$ is a subgraph of $\Gamma'$ and such that $V\Gamma' = V\Gamma$*

*Proof.* Let $x \in X$ be a letter. Denote by $i_x$ the number of vertices in $\Gamma$ which do not have an outgoing edge labelled $x$. Similarly, denote by $j_x$ the number of vertices in $\Gamma$ which do not have an incoming edge labelled $x$. Let $n_x$ be the number of edges in $\Gamma$ labelled $x$. It is clear that $n_x$ is equal to the number of vertices having an incoming edge labeled $x$. Similarly, $n_x$ is equal to the number of vertices having an outgoing edge labeled $x$. Hence

$$\#V\Gamma - n_x = i_x = j_x.$$

Therefore by adding exactly $i_x = j_x$ edges labelled $x$ to $\Gamma$ we can obtained a folded $X$-graph with the same vertex set as $\Gamma$ where each vertex has an incoming edge labelled $x$ and an outgoing edge labelled $x$. Repeating this process for each letter $x \in X$ we obtain an $X$-regular graph $\Gamma'$ with the properties required in Lemma 8.10. □

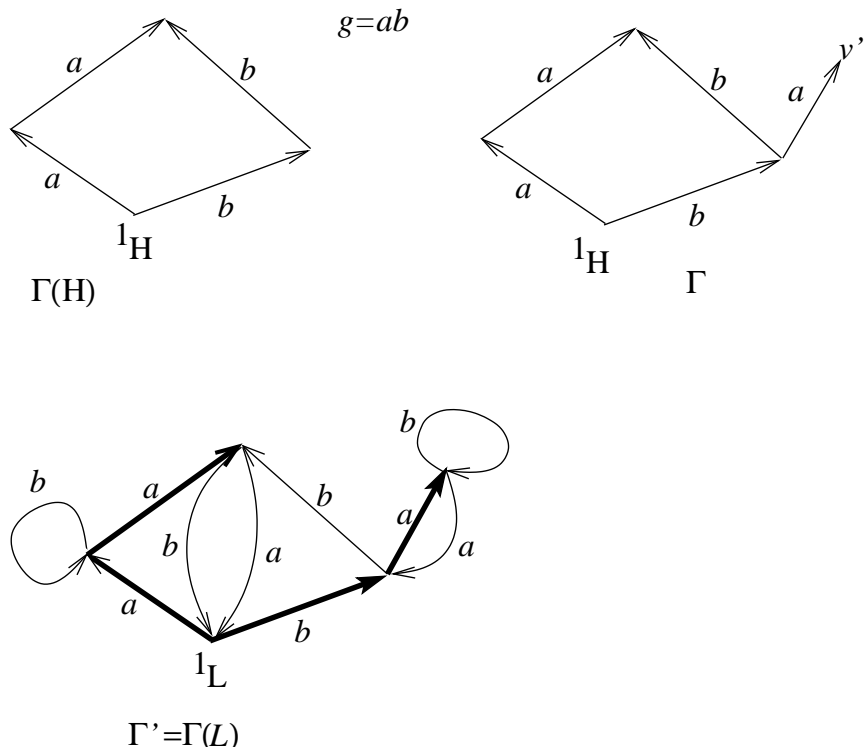

The graph-theoretic criterion of being a subgroup of finite index can also be used to prove the following classical result due to M.Hall [20].

**Theorem 8.11** (Marshall Hall's Theorem). *Let $H$ be a finitely generated subgroup of $F(X)$. Let $g \in F(X)$ be such that $g \notin H$.*

*Then there exists a finitely generated subgroup $K$ of $F(X)$ such that*

1. $L = \langle H, K \rangle = H * K$;
2. *$L$ has finite index in $F(X)$;*
3. $g \notin H$.

*Proof.* Let $H \leq F(X)$ be a finitely generated subgroup. If $H$ has finite index in $F(X)$, then $K = 1$ satisfies the requirement of Theorem 8.11. Suppose now that $H$ has finite index in $F(X)$. Then $\Gamma(H)$ is not $X$-complete.

Write the element $g$ as a freely reduced word in $X$ and construct an $X$-graph $Z$ which is a line-segment subdivided into $|g|$ edges with initial vertex $u$ and terminal vertex $v$ such that the label of $Z$ read from $u$ to $v$ is precisely $g$.

We join the graphs $\Gamma(H)$ and $Z$ by amalgamating $1_H$ and $u$ into a single vertex. Let $\Gamma_0$ be the resulting $X$-graph. We then perform all possible $X$-foldings on $\Gamma_0$ and denote the resulting graph by $\Gamma$. Note that $\Gamma(H)$ is canonically a subgraph of $\Gamma$. We can also think about $\Gamma$ as the result of "wrapping" onto $\Gamma(H)$ the maximal initial segment of $g$ which can be red in $\Gamma(H)$ starting at $1_H$. Let $v'$ be the image of $v$ in $\Gamma$. Note that $v' \neq 1_H$ since $g \notin H$. The graph $\Gamma$ is finite, connected and $X$-folded by construction.

By Lemma 8.10 there exists a connected $X$-regular graph $\Gamma'$ containing $\Gamma$ as a subgraph and having the same vertex set as $\Gamma$.

Thus $L = L(\Gamma', 1_H)$ is a subgroup of finite index (equal to $\#V\Gamma$) in $F(X)$. Since $\Gamma(H)$ is a subgraph of $\Gamma$ and hence of $\Gamma'$, Proposition 6.2 implies that $H$ is a free factor of $L$. Moreover, the path labelled $g$ with origin $1_H$ in $\Gamma'$ terminates at the vertex $v' \neq 1_H$. Hence $g \notin L$, as required. Theorem 8.11 is proved. $\qquad \blacksquare$

It is easy to see that our proof of M.Hall's Theorem is algorithmically effective. An example of how such an algorithm works is given in Figure 7.

**Remark 8.12.** Theorem 8.11 implies that a free group $F(X)$ of finite rank is *subgroup separable*, that is to say any finitely generated subgroup $H$ of $F(X)$ is equal to the intersection of finite index subgroups of $F(X)$ containing $H$. This is an important and nontrivial property of free groups.

**Corollary 8.13.** *There exists an algorithm which, given a finite set of freely reduced words $h_1, \ldots, h_m \in F(X)$ produces the following for the subgroup $H = \langle h_1, \ldots, h_m \rangle \leq F(X)$:*

1. *The graph of a finitely generated subgroup $L \leq F(X)$ such that $H$ is a free factor of $L$ and $L$ is a subgroup of finite index in $F(X)$.*
2. *A free basis $Y$ of $L$ decomposed as a disjoint union $Y = Y_H \cup Y_C$, where $Y_H$ is a free basis of $H$. Thus $L = F(Y) = H * C$ where $C = F(Y_C)$.*
3. *The index of $L$ in $F(X)$.*

*Proof.* We first construct the graph $\Gamma(H)$ and then complete it to an $X$-regular folded graph $\Delta$ without adding any new vertices. Thus $(\Delta, 1_H) = (\Gamma(L), 1_L)$ and $L$ has finite index in $F(X)$. In fact this index is equal to the number of vertices in $\Gamma(H)$.

We then find a spanning tree $T$ in $\Gamma(H)$. Since the vertex sets of $\Gamma(H), \Gamma(L)$ are the same, $T$ is also a spanning tree in $\Gamma(L)$. We then make a free basis $Y = Y_T$ for $L$ as described in Lemma 6.1. The set $Y$ naturally decomposes into a disjoint union $Y = Y_H \cup Y_C$, where $Y_H$ is a free basis of $H$ (those elements $[e]$ of $Y_T$ corresponding to $e \in E^+(\Gamma(H) - T)$ form $Y_H$). Thus $L = H * C$ (where $C = F(Y_C)$) as required. $\qquad \blacksquare$

It turns out that $X$-regular graphs are essential for describing normal subgroups of $F(X)$.

**Theorem 8.14** (Normal subgroups). *Let $H \leq F(X)$ be a nontrivial subgroup of $F(X)$. Then $H$ is normal in $F(X)$ if and only if the following conditions are satisfied:*

1. *The graph $\Gamma(H)$ is $X$-regular (so that there are no degree-one vertices in $\Gamma(H)$ and hence $\Gamma(H) = Type(\Gamma(H))$).*
2. *For any vertex $v$ of $\Gamma(H)$ the based $X$-digraphs $(\Gamma(H), 1_H)$ and $(\Gamma(H), v)$ are isomorphic (that is $L(\Gamma(H), v) = H$).*

*Proof.* Suppose first that conditions (1) and (2) are satisfied. Let $g \in F(X)$ be an arbitrary freely reduced word. Since $\Gamma(H)$ is folded and $X$-regular, there exists a unique path $p$ in $\Gamma(H)$ with origin $1_H$ and label $g$. Let $v$ be the terminal vertex of $p$. By Lemma 7.5 $L(\Gamma(H), v) = gHg^{-1}$. On the other hand by assumption (2) $L(\Gamma(H), v) = H$. Thus $gHg^{-1} = H$. Since $g \in F(X)$ was chosen arbitrarily, this implies that $H$ is normal in $F(X)$.

Suppose now that $H \neq 1$ is normal in $F(X)$. Note that the degree of $1_H$ in $\Gamma(H)$ is at least two. Indeed, suppose there is only one edge with origin $1_H$ and let $x$ be the label of this edge. Then any element of $H$ (considered as a freely reduced word) has the form $xwx^{-1}$. Let $h = xwx^{-1} \in H$ be one such nontrivial element. The cyclically reduced form $q$ of $h$ is conjugate to $h$ and hence belongs to $H$ (since $H$ is normal). However, a cyclically reduced word cannot have the form $x \ldots x^{-1}$.

Thus the degree of $1_H$ is at least two and so $Type(\Gamma(H)) = \Gamma(H)$.

Let $g \in F(X)$ be an arbitrary freely reduced word. We claim that there is a path $p_g$ in $\Gamma(H)$ with origin $1_H$ and label $g$. If this is not the case then by Lemma 7.6 the graph of the subgroup $gHg^{-1}$ has $1_{gHg^{-1}}$ as a degree-one vertex. However, $H = gHg^{-1}$ since $H$ is normal, yielding a contradiction with the above observations about $\Gamma(H)$.

Thus such a path $p_g$ indeed exists. Since $g \in F(X)$ was chosen arbitrarily and the graph $\Gamma(H)$ is folded and connected, this implies that $\Gamma(H)$ is $X$-regular. Moreover, if $v_g$ is the terminal vertex of $p_g$ then by Lemma 7.5 $L(\Gamma(H), v_g) = gHg^{-1} = H$, so that $(\Gamma(H), v_g) \cong (\Gamma(H), 1_H)$. It is clear that any vertex $v$ of $\Gamma(H)$ can be obtained as $v_g$ for some $g \in F(X)$ and so all the conditions of Theorem 8.14 are satisfied. □

**Corollary 8.15.** *There exists an algorithm which, given finitely many freely reduced words $h_1, \ldots, h_m \in F(X)$ decides whether or not $H = \langle h_1, \ldots, h_m \rangle$ is a normal subgroup of finite index in $F(X)$.*

*Proof.* We first construct the graph $\Gamma(H)$ and check whether it is $X$-regular. If not, then the index of $H$ is infinite. If yes, we list the vertices of $\Gamma(H)$ and for each vertex $v$ check if the graphs $(\Gamma(H), 1_H)$ and $(\Gamma(H), v)$ are isomorphic. □

## 9. Intersections of subgroups

One of the most interesting applications of the folded graphs technique is for computing the intersection of two subgroups of a free group.

**Definition 9.1** (Product-graph). Let $\Gamma$ and $\Delta$ be $X$-digraphs. We define the *product-graph* $\Gamma \times \Delta$ as follows.

The vertex set of $\Gamma \times \Delta$ is the set $V\Gamma \times V\Delta$.

For a pair of vertices $(v, u), (v', u') \in V(\Gamma \times \Delta)$ (so that $v, v' \in V\Gamma$ and $u, u' \in V\Delta$ and a letter $x \in X$ we introduce an edge labeled $x$ with origin $(v, u)$ and terminus $(v', u')$ provided there is an edge, labeled $x$, from $v$ to $v'$ in $\Gamma$ and there is an edge, labeled $x$, from $u$ to $u'$ in $\Delta$.

Thus $\Gamma \times \Delta$ is an $X$-digraph. In this situation we will sometimes denote a vertex $(v, u)$ of $\Gamma \times \Delta$ by $v \times u$.

The following trivial observation follows immediately from the definition.

**Lemma 9.2.** *Suppose $\Gamma$ and $\Delta$ are folded $X$-digraphs. Then $\Gamma \times \Delta$ is also a folded $X$-digraph. Moreover, if $v \times u$ is a vertex of $\Gamma \times \Delta$, then the degree of $v \times u$ is less than or equal to the minimum of the degree of $v$ in $\Gamma$ and the degree of $u$ in $\Delta$.*

Next we establish an important connection between product-graphs and intersections of subgroups.

**Lemma 9.3.** *Let $\Gamma$ and $\Delta$ be folded $X$-digraphs. Let $H = L(\Gamma, v)$ and $K = L(\Delta, u)$ for some vertices $v \in V\Gamma$ and $u \in V\Delta$. Let $y = (v, u) \in V(\Gamma \times \Delta)$.*

*Then $\Gamma \times \Delta$ is folded and $L(\Gamma \times \Delta, y) = H \cap K$.*

*Proof.* Let $w$ be a freely reduced word in $X$.

By definition of $\Gamma \times \Delta$ there exists a path in $\Gamma \times \Delta$ from $y$ to $y$ with label $w$ if and only if there is a reduced path in $\Gamma$ from $v$ to $v$ with label $w$ and there is a reduced path in $\Delta$ from $u$ to $u$ with label $w$. This implies $L(\Gamma \times \Delta, y) = L(\Gamma, v) \cap L(\Delta, u) = H \cap K$, as required. □

**Proposition 9.4.** *Let $H \leq F(X)$ and $K \leq F(X)$ be two subgroups of $F(X)$. Let $G = H \cap K$. Let $\Gamma(H) \times_c \Gamma(K)$ be the connected component of $\Gamma(H) \times \Gamma(K)$ containing $1_H \times 1_K$. Let $\Delta = Core(\Gamma(H) \times \Gamma(K), 1_H \times 1_K)$.*

*Then $(\Gamma(G), 1_G) = (\Delta, 1_H \times 1_K)$.*

*Proof.* By Lemma 9.3 the graph $\Gamma(H) \times_c \Gamma(K)$ is folded and $L(\Gamma(H) \times_c \Gamma(K), 1_H \times 1_K) = H \cap K$. Taking the core does not change the language of a graph, so $L(\Delta, 1_H \times 1_K) = H \cap K = G$.

Obviously, $\Delta$ is a folded connected graph which is a core graph with respect to $1_H \times 1_K$. Hence $(\Delta, 1_H \times 1_K) = (\Gamma(G), 1_G)$. □

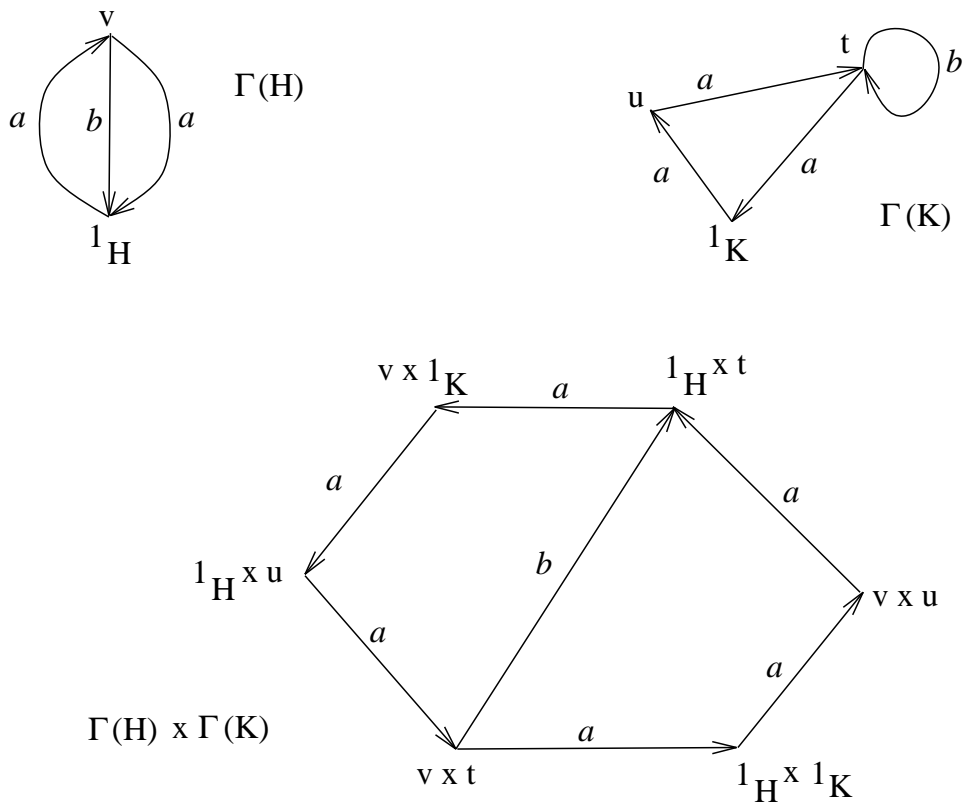

FIGURE 8. Product of two graphs. Here we compute $\Gamma(H) \times \Gamma(K)$ where $H = \langle ab, b^{-1}a \rangle \leq F(a,b)$ and $K = \langle a^3, a^{-1}ba \rangle \leq F(a,b)$. In this example the product is connected and is a core-graph, although in general this need not be the case.

**Corollary 9.5.** *There exists an algorithm which, given finitely many freely reduced words*

$$h_1, \ldots, h_s, k_1, \ldots, k_m \in F(X),$$

*finds the rank and a Nielsen-reduced free basis of the subgroup $\langle h_1, \ldots, h_s \rangle \cap \langle k_1, \ldots, k_m \rangle$ of $F(X)$. In particular, this algorithm determines whether or not $\langle h_1, \ldots, h_s \rangle \cap \langle k_1, \ldots, k_m \rangle = 1$.*

**Corollary 9.6** (Howson Property). [22] *The intersection of any two finitely generated subgroups of $F(X)$ is again finitely generated.*

*Proof.* Let $H, K$ be finitely generated subgroups of $F(X)$ and let $G = H \cap K$. By Lemma 5.4 the graphs $\Gamma(H)$ and $\Gamma(K)$ are finite. Therefore the graph $\Gamma(H) \times \Gamma(K)$ is finite. By Lemma 9.3 $G = L(\Gamma(H) \times \Gamma(K), 1_H \times 1_K)$. Therefore by Lemma 6.1 the group $G$ is finitely generated. $\square$

It turns out that even those components of $\Gamma(H) \times \Gamma(K)$ which don't contain $1_H \times 1_K$, carry some interesting information about $H$ and $K$.

**Proposition 9.7.** *Let $H, K \leq F(X)$. Let $g \in F(X)$ be such that the double cosets $KgH$ and $KH$ are distinct. Suppose that $gHg^{-1} \cap K \neq 1$.*

*Then there is a vertex $v \times u$ in $\Gamma(H) \times \Gamma(K)$ which does not belong to the connected component of $1_H \times 1_K$ such that the subgroup $L(\Gamma(H) \times \Gamma(K), v \times u)$ is conjugate to $gHg^{-1} \cap K$ in $F(X)$.*

*Proof.* Let $g = yz$ where $z$ is the largest terminal segment of $g$ such that $z^{-1}$ is the label of a path in $\Gamma(H)$ with origin $1_H$. Denote this path by $\sigma$ and the terminal vertex of $\sigma$ by $v$.

If the word $y$ is not the label of a path in $\Gamma(K)$ with origin $1_K$, then by Proposition 7.9 $\langle gHg^{-1}, K \rangle = gHg^{-1} * K$. This implies that $gHg^{-1} \cap K = 1$, contrary to our assumptions.

Thus $y$ is the label of a path $\tau$ in $\Gamma(K)$ from $1_K$ to some vertex $u$.

By Lemma 7.5 $L(\Gamma(H), v) = zHz^{-1}$ and $L(\Gamma(K), u) = y^{-1}Ky$. Note also that

$$y(zHz^{-1} \cap y^{-1}Ky)y^{-1} = yzHz^{-1}y^{-1} \cap K = gHg^{-1} \cap K$$

and therefore the subgroups $gHg^{-1} \cap K$ and $zHz^{-1} \cap y^{-1}Ky$ are conjugate in $F(X)$.

By Lemma 9.3 we have

$$L(\Gamma(H) \times \Gamma(K), v \times u) = zHz^{-1} \cap y^{-1}Ky \cong gHg^{-1} \cap K,$$

as required.

It remains to show that $v \times u$ does not belong to the connected component of $1_H \times 1_K$ in $\Gamma(H) \times \Gamma(K)$. Suppose this is not the case.

Then there exist a reduced path $p_v$ in $\Gamma(H)$ from $1_H$ to $v$ and a reduced path $p_u$ in $\Gamma(K)$ from $1_K$ to $u$ such that their labels are the same, that is $\mu(p_v) = \mu(p_u) = \alpha$.

Note that $p_v\sigma^{-1}$ is a path in $\Gamma(H)$ from $1_H$ to $1_H$ and therefore $\overline{\mu(p_v\sigma^{-1})} = \overline{\alpha z} = \alpha \cdot z = h \in H$. Similarly, $\tau p_u^{-1}$ is a path in $\Gamma(K)$ from $1_K$ to $1_K$ and hence $y \cdot \alpha^{-1} = k \in K$. Thus we have

$$g = y \cdot z = y \cdot \alpha^{-1} \cdot \alpha \cdot z = k \cdot h \in KH,$$

contrary to our assumption that $KgH \neq KH$. $\qquad\blacksquare$

The converse of the above statement also holds:

**Proposition 9.8.** *Let $H, K \leq F(X)$ be two subgroups of $F(X)$.*

*Then for any vertex $v \times u$ of $\Gamma(H) \times \Gamma(K)$ the subgroup $L(\Gamma(H) \times \Gamma(K), v \times u)$ is conjugate to a subgroup of the form $gHg^{-1} \cap K$ for some $g \in F(X)$.*

*Moreover, if $v \times u$ does not belong to the connected component of $1_H \times 1_K$, then the element $g$ can be chosen so that $KgH \neq KH$.*

*Proof.* Let $p_v$ be a reduced path in $\Gamma(H)$ from $1_H$ to $v$ with label $\sigma$. Similarly, let $p_u$ be a reduced path in $\Gamma(K)$ from $1_K$ to $u$ with label $\tau$. As we have shown in Lemma 7.5, $L(\Gamma(H), v) = \sigma^{-1}H\sigma$ and $L(\Gamma(K), u) = \tau^{-1}K\tau$. Therefore by Lemma 9.3

$$L(\Gamma(H) \times \Gamma(K), v \times u) = \sigma^{-1}H\sigma \cap \tau^{-1}K\tau \text{ conjugate to } \tau\sigma^{-1}H\sigma\tau^{-1} \cap K,$$

and $g = \tau\sigma^{-1}$ satisfies the requirement of the proposition.

Suppose now that $v \times u$ does not belong to the connected component of $1_H \times 1_K$ in $\Gamma(H) \times \Gamma(K)$ but $g = \tau \cdot \sigma^{-1} \in KH$. Thus $\tau \cdot \sigma^{-1} = kh$ for some $k \in K$, $h \in H$ and therefore

$$k^{-1} \cdot \tau = h \cdot \sigma.$$

Let $\alpha$ be the freely reduced form of the element $k^{-1} \cdot \tau = h \cdot \sigma$. Recall that $k^{-1} \in K$ and so $k^{-1}$ is the label of a reduced $p_1$ in $\Gamma(K)$ from $1_K$ to $1_K$. Then $p_1p_u$ is a path in $\Gamma(K)$ whose label freely reduces to $\alpha$. Therefore there is a reduced path $p_1'$ in $\Gamma(K)$ from $1_K$ to $u$ with label $\alpha$. Similarly, since $h \in H$, there is a path $p_2$ in $\Gamma(H)$ from $1_H$ to $1_H$ with label $h$. Hence $p_2p_v$ is a path in $\Gamma(H)$ from $1_H$ to $v$ whose label freely reduces to $\alpha$. Again, it follows that there is a reduced path $p_2'$ in $\Gamma(H)$ from $1_H$ to $v$ with label $\alpha$. Now the definition of $\Gamma(H) \times \Gamma(K)$ implies that there is a path in $\Gamma(H) \times \Gamma(K)$ from $1_H \times 1_K$ to $v \times u$ with label $\alpha$. However, this contradicts our assumption that $v \times u$ does not belong to the connected component of $1_H \times 1_K$.

Thus $g \notin KH$ and $KgH \neq KH$, as required. $\qquad\blacksquare$

The above facts allow us to use the product-graphs to decide whether a finitely generated subgroup of a free group is *malnormal.* .

**Definition 9.9.** Let $H$ be a subgroup of a group $G$. We say that $H$ is a *malnormal* subgroup of $G$ if for any $g \in G - H$

$$gHg^{-1} \cap H = 1.$$

We say that $H$ is *cyclonormal* if for any $g \in G - H$

$$gHg^{-1} \cap H \text{ is cyclic.}$$

(The notion of a cyclonormal subgroup was suggested by D.Wise [49].)

Malnormal subgroups of free groups and word-hyperbolic groups have recently become the object of intensive studies. In particular, malnormality plays an important role in the Combination Theorem for hyperbolic groups [4], [5], [17], [29]. Various examples, where malnormal subgroups play an important role, can be found [28], [2], [23], [24], [26], [27], [47], [48], [49], [50], [8] and other sources.

**Theorem 9.10.** *Let $H \leq F(X)$ be a subgroup. Then $H$ is malnormal in $F(X)$ if and only if every component of $\Gamma(H) \times \Gamma(H)$, which does not contain $1_H \times 1_H$, is a tree.*

*Proof.* Suppose $H$ is malnormal in $G$. Assume that there is a connected component $C$ of $\Gamma(H) \times \Gamma(H)$, which does not contain $1_H \times 1_H$ and which is not a tree. Let $v \times u$ be a vertex of $C$. Since $C$ is not a tree, there exists a nontrivial reduced path from $v \times u$ to $v \times u$ in $C$ and hence $A = L(C, v \times u) \neq 1$. However by Proposition 9.8 there is $g \in F(X)$ such that $HgH \neq HH = H$ (so that $g \notin H$) and such that $gHg^{-1} \cap H$ is conjugate to $A$ and thus $gHg^{-1} \cap H \neq 1$. This contradicts our assumption that $H$ is malnormal.

Since for any tree $T$ and any vertex $v$ of $T$ $L(T, v) = 1$, the opposite implication of the theorem is just as obvious and follows immediately from Proposition 9.7. $\square$

**Corollary 9.11.** *There exists an algorithm which, given finitely many freely reduced words $h_1, \ldots, h_s$ in $F(X)$ decides whether or not the subgroup $H = \langle h_1, \ldots, h_s \rangle$ is malnormal in $F(X)$. If $H$ is not malnormal, this algorithm will produce a nontrivial element $g \in F(X) - H$ such that $gHg^{-1} \cap H \neq 1$.*

*Proof.* We first construct the finite graph $\Gamma(H)$ as described in Proposition 7.1 and then build the finite graph $\Gamma(H) \times \Gamma(H)$. It remains to check whether those connected components of $\Gamma(H) \times \Gamma(H)$ which do not contain $1_H \times 1_H$, are trees.

If all these components are trees, $H$ is malnormal in $G$. If some component $\Omega$ not containing $1_H \times 1_H$ is not a tree, then $H$ is not malnormal in $G$. Moreover, in this case we choose a vertex $(v, u) \in \Omega$ and a paths with labels $\sigma, \tau$ in $\Gamma(H)$ from $1_H$ to $v, u$ respectively. Put $g = \tau \sigma^{-1}$. Then by the proof of Proposition 9.8 we have $gHg^{-1} \cap H \neq 1$ and $HgH \neq H1H$, so that $g \notin H$. $\square$

A different proof of the above result was given in [2].

**Definition 9.12** (Immersed subgroup). Let $H = \langle h_1, \ldots, h_s \rangle \leq F(X)$ be a finitely generated subgroup of $F(X)$, where each $h_i$ is a nontrivial freely reduced word.

We say that $H$ is *immersed* in $G$ if for any $1 \leq i, j \leq s$ we have

$$|h_i \cdot h_j|_X = |h_i|_X + |h_j|_X$$

and for any $i \neq j, 1 \leq i, j \leq s$ we have

$$|h_i \cdot h_j^{-1}|_X = |h_i|_X + |h_j|_X$$

**Remark 9.13.** Note that if $H$ is an immersed subgroup then for any freely reduced word $w(y_1, \ldots, y_s)$ in the alphabet $\{y_1^{\pm 1}, \ldots, y_s^{\pm 1}\}$ the word $W$, obtained from $w$ by substituting $y_i$ by $h_i$, is freely reduced in $X \cup X^{-1}$. It is also easy to see that for $i \neq j$ the words $h_i$ and $h_j$ have different first letters and different last letters. Moreover, each $h_i$ is cyclically reduced, that is the first letter of $h_i$ is not the inverse of the last letter of $h_i$. For this reason the map $\sigma : \{\pm 1, \ldots, \pm s\} \longrightarrow X \cup X^{-1}$ defined as

$$\sigma : i \mapsto \text{ the first letter of } h_{|i|}^{sign(i)}$$

is injective. Therefore $s \leq \#(X)$.

There are many examples of immersed subgroups (see [25] for details). Recall that $X = \{x_1, \ldots, x_N\}$. For each $i$, $1 \leq i \leq n$ let $h_i$ be a nontrivial freely reduced word in $X$ such that the first and the last letters of $h_i$ are $x_i$. Then the subgroup $H = \langle h_1, \ldots, h_n \rangle$ is obviously immersed in $F(X)$.

The following theorem was first proved, in the context of Stallings' original approach, by D.Wise [50]. An interesting application of this result was given by I.Kapovich [25].

**Theorem 9.14.** *Let $H = \langle h_1, \ldots, h_s \rangle \leq F(X)$ be an immersed subgroup of $F(X)$. Then $H$ is cyclonormal in $F(X)$.*

*Proof.* Consider the graph $\Gamma$ which is a wedge of $s$ circles with labels $h_1, \ldots, h_s$ joined at a single vertex (which we denote $1_H$). The fact that $H$ is immersed in $F(X)$ implies (and, moreover, is equivalent to the condition) that $\Gamma$ is already a folded $X$-digraph. Hence $\Gamma = \Gamma(H)$.

We know from Proposition 9.7 that for any $g \notin H$ the subgroup $gHg^{-1} \cap H$ is either trivial or isomorphic to $L(C, v \times u)$ for some connected component $C$ of $\Gamma(H) \times \Gamma(H)$ which does not contain $1_H \times 1_H$.

Thus to prove cyclonormality of $H$ it is enough to establish that for each such component $C$ and any vertex $v \times u$ of $C$ the subgroup $L(C, v \times u)$ is cyclic (note trivial group is cyclic).

Note that by Lemma 9.2 the degree of a vertex $a \times b$ in a product-graph $A \times B$ is bounded by the minimum of the degrees of $a$ and $b$ in $A$ and $B$ respectively. As we observed before, the degree of every vertex, other than $1_H$ in $\Gamma(H)$ is equal to two. Hence for each vertex $v \times u \neq 1_H \times 1_H$ in $\Gamma(H) \times \Gamma(H)$ the degree of $v \times u$ is at most two. Thus $C$ is a finite connected graph where the degree of each vertex is at most two. Therefore $C$ is either a segment or a circle. In the first case the language of $C$ (with respect to any of its vertices) is trivial and in the second case this language is an infinite cyclic subgroup of $F(X)$. This implies that $H$ is cyclonormal in $F(X)$. $\qquad\square$

**Corollary 9.15.** *There is an algorithm which, given finitely many nontrivial freely reduced words*

$$h_1, \ldots, h_s \in F(X),$$

*decides whether or not the subgroup $H = \langle h_1, \ldots, h_s \rangle \leq F(X)$ is cyclonormal.*

*Proof.* If follows from Proposition 9.7 that if $g \notin H$ then $gHg^{-1} \cap H$ is either trivial or isomorphic to $L(C, v \times u)$ where $C$ is a connected component of $\Gamma(H) \times \Gamma(H)$ not containing $1_H \times 1_H$ and where $v \times u$ is a vertex of $C$. Recall also that the isomorphism type of $L(C, v \times u)$ depends only on $C$ and not on the choice of $v \times u$ in $C$.

Thus to check whether $H$ is cyclonormal we construct $\Gamma(H)$, $\Gamma(H) \times \Gamma(H)$ and for every connected component $C$ of $\Gamma(H) \times \Gamma(H)$, not containing $1_H \times 1_H$ we find a spanning tree $T_C$ in $C$ and count the number $n_C$ of positive edges in $C - T_C$ (this number is equal to the rank of the free group $L(C, v \times u)$). If for each $C$ $n_C \leq 1$, the subgroup $H$ is cyclonormal. If for some $C$ we have $n_C \geq 2$, the subgroup $H$ is not cyclonormal. $\qquad\square$

## 10. Principal quotients and subgroups

In this section we want to investigate more carefully the relationship between $\Gamma(H)$ and $\Gamma(K)$ if $K \leq H$. To simplify the matters we will assume that both $H$ and $K$ are finitely generated, although it will be clear that one can consider a more general situation.

**Definition 10.1.** Let $K \leq H \leq F(X)$ be two finitely generated subgroups of $F(X)$. Since $K = L(\Gamma(K), 1_K) \subseteq L(\Gamma(H), 1_H) = H$, for any path $p$ in $\Gamma(K)$ from $1_K$ to $1_K$ there is a unique path $\alpha[p]$ in $\Gamma(H)$ from $1_H$ to $1_H$ such that the label of $\alpha[p]$ is the same as that of $p$.

Define a subgraph $\Gamma_H(K)$ of $\Gamma(H)$ as

$$\Gamma_H(K) = \cup\{\alpha[p] \mid p \text{ is a path in } \Gamma(K) \text{ from } 1_K \text{ to } 1_K\}.$$

In fact, constructing $\Gamma_H(K)$ can be considerably simplified, as the following lemma shows.

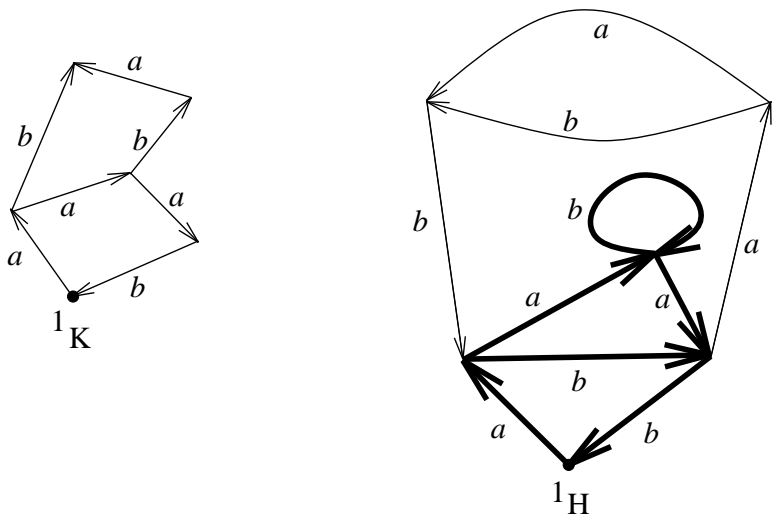

FIGURE 9. An example of $\Gamma_H(K)$. Here $K \leq H \leq F(a,b)$ and the graph $\Gamma_H(K)$ is highlighted. We suggest that the reader uses Lemma 10.2 to verify the validity of this picture

**Lemma 10.2.** *Let $K \leq H \leq F(X)$ be finitely generated subgroups of $F(X)$. Let $K = \langle k_1, \ldots, k_m \rangle$ where $k_i$ are freely reduced words in $F(X)$. For each $i = 1, \ldots, m$ let $p_i$ be the loop at $1_K$ in $\Gamma(K)$ with label $k_i$. Then*

$$\Gamma_H(K) = \cup_{i=1}^n \alpha[p_i]$$

The above lemma shows that in practical terms it is very easy to construct $\Gamma_H(K)$, given $\Gamma(K)$ and $\Gamma(H)$. For example, we can choose as a finite generating set for $K$ the set $T_K$ corresponding to some spanning tree in $\Gamma(K)$.

**Lemma 10.3.** *Let $K \leq H \leq F(X)$. Then there is a unique morphism of based $X$-digraphs $\alpha : (\Gamma(K), 1_K) \longrightarrow (\Gamma(H), 1_H)$; for this morphism we have $\alpha(\Gamma(K)) = \Gamma_H(K)$.*

*Proof.* The existence and uniqueness of $\alpha$ follow from Proposition 4.3.

It is obvious from the definition of $\Gamma_H(K)$ and the fact that any edge in $\Gamma(K)$ lies on a reduced path from $1_K$ to $1_K$ that $\alpha(\Gamma(K)) = \Gamma_H(K)$.     □

**Lemma 10.4.** *Let $K \leq H \leq F(X)$. Then the graph $\Gamma_H(K)$ is a connected folded graph which is a core graph with respect to $1_H$.*

*Proof.* It is obvious that $\Gamma_H(K)$ is folded, connected and that it contains $1_H$. Let $\alpha : \Gamma(K) \longrightarrow \Gamma(H)$ be the unique morphism such that $\alpha(1_K) = 1_H$, so that $\alpha(\Gamma(K)) = \Gamma_H(K)$ (the existence of $\alpha$ follows from Lemma 10.3).

Assume that $\Gamma_H(K)$ is not a core-graph with respect to $1_H$ and that there is a vertex $u \neq 1_H$ of degree one on $\Gamma_H(K)$. Let $e'$ be the unique edge of $\Gamma_H(K)$ with terminus $u$. Let $e$ be an edge of $\Gamma(K)$ such that $\alpha(e) = e'$. Since $\Gamma(K)$ is a core graph with respect to $1_K$, there exists a reduced loop $p$ at $1_K$ passing through $p$. Denote the label of $p$ by $w$. Since $p$ is reduced and $\Gamma(K)$ is folded, $w$ is freely reduced. The path $\alpha(p)$ is a loop at $1_H$ in $\Gamma_H(K)$ with label $w$ and $\alpha(p)$ passes through $e'$. Since $w$ is freely reduced and $\Gamma_H(K)$ is folded, the path $\alpha(p)$ is reduced. However, this contradicts our assumption that $u$ has degree one in $\Gamma_H(K)$.     □

Recall that by Lemma 4.4 (morphism factorization) an epimorphic image of an $X$ digraph $\Gamma$ can be obtained by first identifying some vertex subsets of $\Gamma$ into single vertices and then folding some edges with the same initial vertices and the same terminal vertices (so that the vertex set of the graph does not change at this step).

**Definition 10.5** (Principal quotients). Let $K$ be a finitely generated subgroup of $F(X)$. Since $\Gamma(K)$ is finite, there exist only finitely many based $X$-digraphs $(\Gamma_1, v_1), \dots (\Gamma_s, v_s)$ such that:
(a) each $\Gamma_i$ is a finite connected folded $X$-digraph which is a core-graph with respect to $v_i$;
(b) for each $i = 1, \dots, s$ there is an epimorphism $f : \Gamma(K) \longrightarrow \Gamma_i$ with $f(1_K) = v_i$.

These graphs $(\Gamma_1, v_1), \dots (\Gamma_s, v_s)$ are called *principal quotients* of $\Gamma(K)$.

**Remark 10.6.** Note that if $(\Gamma, v)$ is a principal quotient of $\Gamma(K)$ then $(\Gamma, v) = (\Gamma(H), 1_H)$ for some finitely generated subgroup $H$ of $F(X)$ with $K \leq H$. Moreover, in this case $\Gamma_H(K) = \Gamma$.

**Lemma 10.7.** *Let $F$ be a free group (of possibly infinite rank) and let $K$ be a finitely generated free factor of $F$. Let $H \leq F$ be a finitely generated subgroup such that $K \leq H$. Then $K$ is a free factor of $H$.*

*Proof.* Let $Y = Y_1 \sqcup Y_2$ be a free basis of $F$ such that $Y_1$ is a finite free basis of $K$. Since $H$ is finitely generated and involves only finitely many letters from $Y_2$, we may assume that $Y_2$ is also finite. It is clear that in this case the graph $\Gamma^Y(K)$ is just a wedge of circles, labeled by elements of $Y_1$, wedged at a vertex $1_K$. This graph clearly has no epimorphic images other than itself, so $\Gamma^Y(K) = \Gamma_H^Y(K) \subseteq \Gamma^Y(H)$. Since $(\Gamma^Y(K), 1_K)$ is a subgraph of $(\Gamma^Y(H), 1_H)$, this implies that $K$ is a free factor of $H$, as required. $\qquad\square$

## 11. Algebraic extensions

In this section we will develop an elementary analogue of the field extension theory for subgroups of free groups.

**Definition 11.1** (Free and algebraic extensions). Let $K \leq H \leq F(X)$. In this case we will say that $H$ is an *extension* of $K$. We say that $H$ is a *free extension* of $K$ if there is a subgroup $C \leq F(X), C \neq 1$ and a subgroup $K'$ of $F(X)$ such that $K \leq K'$ and $H = K * C$.

If an extension $K \leq H$ is not free, we call it *algebraic*.

This definition is motivated by the analogy with field extensions. Indeed, if $K \leq K'$ and $H = K' * C$ then $H$ is obtained by adding to $K'$ several "purely transcendental" elements, namely a free basis of $C$. That is why it is reasonable to call such extensions $K \leq H$ "free".

**Theorem 11.2.**  1. *Let $K \leq F(X)$ be a finitely generated subgroup. If $K \leq H$ (where $H \leq F(X)$) is an algebraic extension then $(\Gamma(H), 1_H)$ is a principal quotient of $\Gamma(K)$. In particular, $H$ is finitely generated.*

2. *For a finitely generated $K \leq F(X)$ there are only finitely many subgroups $H \leq F(X)$ such that $K \leq H$ is an algebraic extension.*

3. *Let $K \leq H$ be finitely generated subgroups of $F(X)$ and suppose that the extension $K \leq H$ is free. Then there exist subgroups $K', C$ such that $K \leq K'$, $H = K' * C$ and $(\Gamma(K'), 1_{K'})$ is a principal quotient of $\Gamma(K)$.*

4. *Let $K \leq H$ be finitely generated subgroups of $F(X)$. Then there is a free factor $K'$ of $H$ such that $K \leq K'$ is an algebraic extension.*

*Proof.* (1) Let $K \leq H$ (where $H \leq F(X)$) be a algebraic extension. Then there is a canonical morphism $\alpha : \Gamma(K) \longrightarrow \Gamma(H)$ such that $\alpha(1_K) = 1_H$ and $\alpha(\Gamma(K)) = \Gamma_H(K)$. Since $K$ is finitely generated, the graphs $\Gamma(K)$ and $\Gamma_H(K)$ are finite. Denote $K' = L(\Gamma_H(K), 1_H)$. Then $K'$ is a finitely generated subgroup containing $K$. Moreover, since $\Gamma_H(K)$ is a subgraph of $\Gamma(H)$, the subgroup $K'$ is a free factor of $H$. That is, $H = K' * C$. Since the extension $K \leq H$ is algebraic, $C = 1$ and so $H = K'$. Since $\Gamma_H(K)$ is a connected subgraph of $\Gamma(H)$ which is a core graph with respect to $1_H$ and such that $L(\Gamma_H(K), 1_H) = H$, we have $\Gamma_H(K) = \Gamma(H)$. Thus $(\Gamma(H), 1_H)$ is a principal quotient of $\Gamma(K)$.

(2) If $K$ is a finitely generated subgroup of $F(X)$ then $\Gamma(K)$ is finite and hence has only finitely many principal quotients. Therefore by (1) there are only finitely many algebraic extensions of $K$ in $F(X)$.

(3) Let $K \leq H$ be finitely generated subgroups of $F(X)$. Let $K', C$ be such that $H = K' * C$, $K \leq K'$ and $C$ is of maximal possible rank.

We claim that $\Gamma(K', 1_{K'})$ is a principal quotient of $\Gamma(K)$. Indeed, since $K \leq K'$, we can find the subgraph $\Gamma_{K'}(K)$ of $\Gamma(K')$, which is a principal quotient of $\Gamma(K)$. If $\Gamma_{K'}(K)$ is a proper subgraph in $\Gamma(K')$ then by Proposition 6.2 the subgroup $D = L(\Gamma_{K'}(K), 1_{K'})$ is a proper free factor of $K'$ and $K' = D * D'$ with $D' \neq 1$. On the other hand $K \leq D$ by construction. Thus $K \leq D$ and $H = K' * C = (D * D') * C = D * (D' * C)$. Since the rank of $D' * C$ is greater than the rank of $C$, we obtain a contradiction with the choice of $C$.

Thus $\Gamma_{K'}(K) = \Gamma(K')$, and therefore $(\Gamma(K'), 1'_K)$ is indeed a principal quotient of $\Gamma(K)$.

(4) Let $K \leq H$ be finitely generated subgroups of $F(X)$.

Choose a free factor $K'$ of $H$ of the smallest possible rank such that $K'$ contains $K$. Thus $H = K' * C$ for some (possibly trivial) $C \leq H$.

We claim that $K \leq K'$ is an algebraic extension. Indeed, suppose not. Then there exist nontrivial groups $K'', C'$ such that $K \leq K''$ and $K' = K'' * C'$. This means that $H = K'' * (C' * C)$ and $K \leq K''$. Moreover $C' \neq 1$ and $K' = K'' * C'$ imply that $rk(K'') < rk(K')$. Thus $K''$ is a free factor of $H$ containing $K$ and the rank of $K''$ is smaller than the rank of $K'$. This contradicts the choice of $K'$. $\qquad\square$

**Theorem 11.3** (Decidability of extension type).    1. *There is an algorithm which, given finitely many freely reduced words $k_1, \ldots, k_m$ in $F(X)$ finds all principal quotients of $\Gamma(K)$, where $K = \langle k_1, \ldots, k_m \rangle$.*

  2. *There is an algorithm which, given finitely many freely reduced words $k_1, \ldots, k_m, h_1, \ldots, h_s$ in $F(X)$ decides whether or not the subgroup $H = \langle h_1, \ldots, h_s \rangle$ is an algebraic extension of the subgroup $K = \langle k_1, \ldots, k_m \rangle$.*

  3. *There is an algorithm which, given finitely many elements $f_1, \ldots, f_s \in F(X)$, finds all possible algebraic extensions of the subgroup $\langle f_1, \ldots, f_s \rangle$ in $F(X)$.*

*Proof.* (1) First we construct the graph $\Gamma(K)$, as described in Proposition 7.1. We compute all the principal quotients $(\Gamma_1, v_1), \ldots, (\Gamma_r, v_r)$ of $\Gamma(K)$ as follows. We list all possible partitions of the vertex set of $\Gamma(K)$ into finitely many disjoint subsets. For each partition we identify each of this subsets into a single vertex. We then perform all possible foldings on the resulting graph which do not change the number of vertices. If the final graph is a folded graph which is a core graph with respect to the image of $1_K$, it is a principal quotient of $\Gamma(K)$. Lemma 4.4 ensures that all principal quotients of $\Gamma(K)$ can be obtained in such a way.

(2) First we construct $\Gamma(H)$ and $\Gamma(K)$ and check whether $K$ is contained in $H$. If not, then $H$ is not an extension of $K$ and thus not an algebraic extension of $K$.

Suppose $K \leq H$. We then find all the principal quotients $(\Gamma_1, v_1), \ldots, (\Gamma_r, v_r)$ of $\Gamma(K)$, as described in part (1). After that, for each $i = 1, \ldots, r$ we apply the algorithm of Whitehead (see [32], [33]) to check whether $K_i = L(\Gamma_i, v_i)$ is a proper free factor of $H$. If for some $i$ the answer is yes, then obviously $K \leq K_i$, $K_i$ is a proper free factor of $H$ and thus $K \leq H$ is a free extension.

If for each $i = 1, \ldots, r$ the subgroup $K_i$ is not a proper free factor of $H$, then by part (3) of Theorem 11.2 the extension $K \leq H$ is algebraic.

(3) We first construct the graph $\Gamma(K)$ and all of its principal quotients $(\Gamma_1, v_1), \ldots, (\Gamma_r, v_r)$. Denote $H_i = L(\Gamma_i, v_i)$ for $i = 1, \ldots, r$, so that $K \leq H_i$ for each $i$. By part (1) of Theorem 11.2 if $K \leq H$ is an algebraic extension then $(\Gamma(H), 1_H)$ is a principal quotient of $\Gamma(K)$.

Thus it suffices to check for each $i = 1, \ldots, r$ whether $K \leq H_i$ is algebraic. The set of those $H_i$, for which the answer is "yes", gives us all algebraic extensions of $K$. $\qquad\square$

**Remark 11.4.** The proof of Theorem 11.3 is one of the few places in this paper where the algorithm we provide is slow. Primarily this is because we need to use the Whitehead algorithm for deciding whether one subgroup is a free factor of another.

**Lemma 11.5.** *Let $K$ be a finitely generated subgroup of $F(X)$. Suppose $X = X_1 \cup X_2$ and suppose that $K \leq F(X_1)$ (that is the graph $\Gamma(K)$ does not have any edges labeled by elements of $X_2^{\pm 1}$). Let $H \leq F(X)$ be a subgroup containing $K$ such that the extension $K \leq H$ is algebraic.*

*Then $H \leq F(X_1)$ (so that the graph $\Gamma(H)$ does not have any edges labeled by elements of $X_2^{\pm 1}$).*

**Theorem 11.6.** *Let $K$ be a finitely generated subgroup of $F(X)$.*

1. *If $K \leq H \leq F(X)$ and $|H : K| < \infty$ then $K \leq H$ is an algebraic extension.*
2. *If $K \leq H$ and $H \leq Q$ are algebraic then $K \leq Q$ is algebraic.*
3. *If $K \leq H_1$ and $K \leq H_2$ are algebraic then $K \leq H$ is also algebraic, where $H = \langle H_1, H_2 \rangle$.*

*Proof.* (1) This statement follows directly from Definition 11.1.

(2) Suppose the statement fails and $K \leq Q$ is a free extension. Note that since $K$ is finitely generated, $K \leq H$ is algebraic and $H \leq Q$ is algebraic, both $H$ and $Q$ are also finitely generated. Let $Q = K' * C$ where $K \leq K'$ and $C \neq 1$. Let $Y'$ be a free basis of $K'$ and let $Y$ be a free basis of $C$. Then $Z = Y' \cup Y$ is a free basis of $Q$, so that $Q = F(Z)$. Since $K \leq K' = F(Y')$, the graph $\Gamma^Z(K)$ (corresponding to $K$ with respect to the basis $Z$ of $Q$) does not have any edges labeled by letters of $Y^{\pm 1}$. Since $K \leq H$ is algebraic, Lemma 11.5 implies that the graph $\Gamma^Z(H)$ also does not have any edges labeled by elements of $Y^{\pm 1}$. This means that $H \leq F(Y')$. Since $Q = F(Y') * F(Y)$ and $Y \neq \emptyset$, the extension $H \leq Q$ is free. However, this contradicts our assumption that $H \leq Q$ is algebraic.

(3) Suppose the statement fails and $K \leq H$ is a free extension. Note that $K, H_1, H_2$ and $H$ are finitely generated. Thus $H = K' * C$ where $K' = F(Y')$ and $C = F(Y)$ with finite $Y, Y'$ and nonempty $Y$. Denote $Z = Y' \cup Y$, so that $H = F(Z)$. Since $K \leq K' = F(Y')$, the graph $\Gamma^Z(K)$ does not have any edges with labels from $Y^{\pm 1}$. Since $K \leq H_1$ and $K \leq H_2$ are algebraic, Lemma 11.5 implies that $H_1 \leq F(Y')$ and $H_2 \leq F(Y')$. Therefore $H = \langle H_1, H_2 \rangle \leq F(Y')$. However, this contradicts our assumption that $H = F(Y') * F(Y)$ with $Y \neq \emptyset$. $\qquad \square$

**Definition 11.7** (Algebraic closure). Let $K \leq F(X)$ be a finitely generated subgroup of $F(X)$. Since there are only finitely many algebraic extensions of $K$ in $F(X)$, Theorem 11.6 implies that there exists the largest algebraic extension $H$ of $K$ (which contains all other such extensions). We call this $H$ *the algebraic closure* of $K$ and denote $H = cl(K)$ (or $H = cl_{F(X)}(K)$). If $K = cl(K)$ we say that $K$ is an *algebraically closed* subgroup of $F(X)$.

**Lemma 11.8.** *There is an algorithm which, given finitely many freely reduced words $h_1, \ldots, h_m \in F(X)$ finds a free basis of the subgroup $cl(H)$ of $F(X)$ where $H = \langle h_1, \ldots, h_m \rangle$.*

*Proof.* This follows directly from decidability of the extension type and the fact that we can effectively find all algebraic extensions of $H$ in $F(X)$ (Theorem 11.3). $\qquad \square$

It turns out that algebraically closed subgroups of $F(X)$ are precisely the free factors of $F(X)$.

**Lemma 11.9.** *Let $K \leq H$ be finitely generated subgroups of $F(X)$ such that $H$ is a free extension of $K$. Let $C \neq 1$ and $K'$ be such that $K \leq K'$, $H = K' * C$ and $C$ is of the largest possible rank. Then $K \leq K'$ is an algebraic extension.*

*Proof.* Suppose the statement of Lemma 11.9 fails and the extension $K \leq K'$ is free. This means that $K \leq K''$ and $K' = K'' * C''$ where $C'' \neq 1$. Therefore

$$H = K' * C = (K'' * C'') * C = K'' * (C'' * C).$$

Since $C'' \neq 1$, $rk(C'' * C) = rk(C'') + rk(C) > rk(C)$, which contradicts the choice of $C$. $\qquad \square$

**Theorem 11.10** (Algebraically closed subgroups). *Let $K$ be a finitely generated subgroup of $F(X)$. Then $K$ is algebraically closed in $F(X)$ if and only if $H$ is a free factor of $F(X)$, that is $F(X) = H * C$ for some $C \leq F(X)$.*

*Proof.* First assume that $K$ is a free factor of $F(X)$. Without loss of generality we may assume that $K = F(X_1)$ where $X = X_1 \cup X_2$. Suppose there $K \leq H \leq F(X)$ where $K \leq H$ is an algebraic extension. Hence by Lemma 11.5 $H \leq F(X_1) = K$, so that $H = K$. Hence $K$ is in fact algebraically closed, as required.

Suppose now that $K$ is algebraically closed and $K \neq F(X)$. Therefore the extension $K \leq F(X)$ is free. Let $K', C$ be such that $K \leq K'$, $F(X) = K' * C$ and $C \neq 1$ is of the biggest rank possible. We claim that $K = K'$. Indeed, if $K \neq K'$ then the extension $K \leq K'$ is free since $K$ is algebraically closed. Thus there is $C' \neq 1$ and $K''$ such that $K \leq K''$ and $K' = K'' * C'$. This means that

$$F(X) = K' * C = (K'' * C') * C = K'' * (C' * C)$$

where $K \leq K''$ and the rank $rk(C' * C) = rk(C') + rk(C)$ is larger than the rank of $C$. This contradicts the choice of $C$. Thus $K = K'$ so that $F(X) = K * C$ and $K$ is a free factor of $F(X)$, as required. □

## 12. A remark on Hanna Neumann's conjecture

By a well-known result of A.Howson (see Corollary 9.6 above) the intersection of any two finitely generated subgroups in a free group is again finitely generated. However, it is interesting to investigate in more detail the connection between the rank of the intersection and the ranks of the two intersecting subgroups. A very important and still open classical conjecture regarding this question was proposed by Hanna Neumann [34]:

**Conjecture 12.1.** *Let $A, B$ be finitely generated subgroups of a free group $F(X)$. Let $C = A \cap B$ and suppose that $C \neq 1$. Then*

$$rk(C) - 1 \leq (rk(A) - 1)(rk(B) - 1).$$

Over the last four decades the Hanna Neumann conjecture has been the subject of extensive research and a great deal is known by now (see for example [12], [35], [10]), although the conjecture itself remains open. The following simple statement says that in order to prove the conjecture it is enough to establish its validity for the case when both $A$ and $B$ are algebraic extensions of $C = A \cap B$.

**Theorem 12.2.** *Suppose the Hanna Neumann conjecture holds for all finitely generated subgroups $A$ and $B$ of a free group $F(X)$ such that $C = A \cap B \neq 1$ and such that both $A$ and $B$ are algebraic extensions of $C = A \cap B$. Then the Hanna Neumann conjecture holds for all finitely generated subgroups $A$ and $B$ of $F(X)$ with nontrivial intersection.*

*Proof.* Suppose the Hanna Neumann conjecture holds for all finitely generated subgroups of $F(X)$ such that both of them are algebraic extensions of their intersection.

Let $A$ and $B$ be two arbitrary finitely generated subgroups of $F(X)$. Let $C = A \cap B$ and suppose that $C \neq 1$. By part (4) of Theorem 11.2 there exist subgroups $A', B'$ such that:

(1) we have $C \leq A' \leq A$, $C \leq B' \leq B$;

(2) the subgroup $A'$ is a free factor of $A$ and the subgroup $B'$ is a free factor of $B$;

(3) both $A'$ and $B'$ are algebraic extensions of $C$.

Note that obviously $C = A' \cap B'$. Observe also that (2) implies that $rk(A') \leq rk(A)$ and $rk(B') \leq rk(B)$. By our assumption the Hanna Neumann conjecture holds for algebraic extensions. Therefore

$$rk(C) - 1 \leq (rk(A') - 1)(rk(B') - 1) \leq (rk(A) - 1)(rk(B) - 1).$$

Since $A$ and $B$ were chosen arbitrarily, this implies the statement of Theorem 12.2. □

## 13. Malnormal closures and isolators

We have already mentioned the importance of malnormal subgroups. Recently the significance of *isolated* subgroups has also been clarified. Recall that a subgroup $H$ of a group $G$ is called *isolated* in $G$ if whenever $g^n \in H, g \in G, n > 0$, we also have $g \in H$. It turns out that isolated subgroups (in particular those of free groups) play a substantial role in the study of equations over free groups and the elementary theory of free groups (see [31], [30]). In this section we obtain some new results regarding malnormal and isolated subgroups of free groups.

Unlike malnormality (see Theorem 9.10), verifying the property of being isolated is not as simple. Nonetheless, it can still be done using the subgroup graph.

**Theorem 13.1.** *There exists an algorithm which, given a finite collection of words $h_1, \ldots, h_l \in F(X)$ decides if the subgroup $H = \langle h_1, \ldots, h_l \rangle$ is isolated in $F(X)$.*

*Proof.* Let $\Gamma = \Gamma(H)$ be the subgroup graph of $H$. Denote $n := \#X$ and $k := \#V\Gamma$.

We will first prove the following

**Claim.** Suppose $g^m \in H$ where $g \in F - H$ and $m > 1$. Then there is $f \in F - H$ such that $f^m \in H$ and $|f| \leq [(2n)^k k^{2m} + 1](k+1) + 2k$.

Let $f$ be the shortest element in $F - H$ such that $f^m \in H$. We write $f$ as $f = ara^{-1}$ where $r$ is cyclically reduced in $X$. Since $f^m = ar^m a^{-1} \in H$, every initial segment of the word $ar^m a^{-1}$ can be read as a label of a path in $\Gamma(H)$ originating at $1_H$. For each $i = 0, 1, \ldots, m$ let $v_i$ be the terminal vertex of a path from $1_H$ in $\Gamma(H)$ labeled $ar^i$. Since $f = ara^{-1} \notin H$, we conclude that $v_0 \neq v_1$. Similarly, since $f^m = ar^m a^{-1} \in H$, we have $v_0 = v_m$. Thus for any word $a'$ which is a label of a path in $\Gamma(H)$ from $1_H$ to $v_0$ we have $a'r(a')^{-1} \notin H$ and $a'r^m(a')^{-1} \in H$. Since $f$ was assumed to be shortest in $F - H$ with $f^m \in H$, this implies that the path from $1_H$ to $v_0$ in $\Gamma(H)$ labeled $a$ has no cycles and hence $|a| \leq k = \#V\Gamma$.

Let $r = b'cb''$ where $c$ is the label of a simple loop in the path labeled $r$ from $v_0$ to $v_1$ in $\Gamma(H)$ (so that $|c| \leq k$). Denote by $u_i$ and $w_i$ the terminal vertices of the paths labeled $b'$ and $b'c$ originating from $v_i$. Thus $c$ is the label of a path from $u_i$ to $w_i$. Moreover, $u_0 = w_0$ by the choice of $c$. Consider the tuple $Q = (c, u_0, u_1, w_1, u_2, w_2, \ldots, u_{m-1}, w_{m-1})$. Since $|c| \leq k$, there are at most $N := (2n)^k k^{2m}$ possibilities for $Q$.

Suppose $|r| \geq (N+1)(k+1)$. Since each path of length $k+1$ in $\Gamma(H)$ contains a cycle, we can write $r$ as $r = b_0 c_0 b_1 c_1 \ldots c_N b_{N+1}$ where each $c_j$ has length at most $k$ and is the label of a nontrivial simple loop in the path labeled $r$ from $v_0$ to $v_1$ in $\Gamma(H)$. Let $Q_j$ be the tuple produced as above for each $j = 0, 1, \ldots, N$. By the choice of $N$ there are some $j < s$ such that $Q_j = Q_s$. Let $z = c_j b_{j+1} \ldots b_s$. Thus $z$ is a subword of $r$ which corresponds to a cycle in each path labeled $r$ from $v_i$ to $v_{i+1}$ for $i = 0, 1, \ldots, m-1$.

Let $r'$ be the word obtained by deleting the subword $z$ from $r$. That is $r' = b_0 c_0 \ldots b_j c_s b_{s+1} \ldots b_{N+1}$. By construction $|r'| < |r|$, $ar'a^{-1} \notin H$ and $a(r')^m a^{-1} \in H$. This contradicts the choice of $f$. Thus $|r| < (N+1)(k+1)$. Since $|a| \leq k$, this implies that $|f| = |ara^{-1}| \leq (N+1)(k+1) + 2k$ and the Claim holds.

Recall that by Proposition 7.10 if $g^m \in H$ then $g^{m_0} \in H$ for some $1 \leq m_0 \leq k = \#V\Gamma$. We can now formulate the algorithm for checking if $H$ is isolated.

First, construct the graph $\Gamma = \Gamma(H)$ and compute $k = \#V\Gamma$. Put $M := [(2n)^k k^{2k} + 1](k+1) + 2k$. Now for each word $f$ in $F(X)$ with $|f| \leq M$ such that $f \notin H$ check if there is $2 \leq m \leq k$ such that $f^m \in H$. If yes then $H$ is not isolated. If no such $f$ exists, then $H$ is isolated. $\qed$

**Remark 13.2.** The above theorem is one of the few places in the present paper where our algorithm has high complexity. In fact, quite a few improvements are possible in the algorithm provided in Theorem 13.1. For example, if a subgroup is malnormal (which is easy to check) then it is also isolated. However, at present we do not know of a fast algorithm (comparable in speed with the malnormality test provided by Theorem 9.10) for determining if a finitely generated subgroup of a free group is isolated.

**Lemma 13.3.** *Let $K$ be a finitely generated subgroup of $F(X)$ which is not malnormal in $F(X)$. Let $g \in F(X)$ be such that $g \notin K$ and $gKg^{-1} \cap K \neq 1$. Then the subgroup $H = \langle K, g \rangle$ is an algebraic extension of $K$. Moreover, the rank of $H$ is less than or equal to the rank of $K$.*

*Proof.* Suppose, on the contrary, that $K \leq H$ is a free extension. Then $K \leq K'$ where $H = K' * C$ and $C \neq 1$. Note that $g^{-1}K'g \cap K' \neq 1$. On the other hand $K'$ is a free factor of $H$ and so malnormal in $H$. Hence $g \in K'$. Since $K \leq K', g \in K'$, we have $H = \langle K, g \rangle \leq K'$. This is impossible since $H = K' * C$ and $C \neq 1$. Thus $K \leq H$ is a algebraic extension of $K$, as required.

Since $H = \langle K, g \rangle$ and $H, K$ are free, then either $rk(H) = rk(K)+1$ and $H = K*\langle g \rangle$, or $rk(H) \leq rk(K)$. The former contradicts our assumption that $g^{-1}Kg \cap K \neq 1$. Thus $rk(H) \leq rk(K)$, as required. $\qed$

**Lemma 13.4.** *Let $K$ be a finitely generated subgroup of $F(X)$. Then there exists a unique subgroup $H$ of $F(X)$ such that $H$ is minimal among algebraic extensions of $K$ which are malnormal in $F(X)$.*

*Proof.* If $K$ is malnormal in $F(X)$, the statement is obvious. Assume now that $K$ is not malnormal.

First we observe that there is at least one malnormal subgroup which is an algebraic extension of $K$.

Let $K = K_1$ and define a sequence $K_1 \leq K_2 \leq \ldots$ as follows. Since $K_1$ is not malnormal, there is $g_2 \in F(X) - K_1$ such that $g_2 K_1 g_2^{-1} \cap K_1 \neq 1$. Put $K_2 = \langle K_1, g_2 \rangle$. Suppose now that $K_1 \leq K_2 \leq \cdots \leq K_i$ are already constructed. If $K_i$ is malnormal in $F(X)$, we terminate the sequence. Otherwise there is $g_{i+1} \in F(X) - K_i$ such that $g_{i+1} K_i g_{i+1}^{-1} \cap K_i \neq 1$. Put $K_{i+1} = \langle K_i, g_{i+1} \rangle$.

**Claim.** The sequence $K_1 \leq K_2 \leq \ldots$ terminates in a finite number of steps. Moreover, each $K_i$ is a algebraic extension of $K$ and $rk(K_i) \leq rk(K)$.

Lemma 13.3 implies that $K_i \leq K_{i+1}$ is a algebraic extension and $rk(K_{i+1}) \leq rk(K_i)$. By transitivity (see Theorem 11.6) this means that $K \leq K_i$ is algebraic for each $i \geq 1$. Since $K_i \neq K_{i+1}$ and there are only finitely many algebraic extensions of $K$, the sequence terminates. Also, obviously $rk(K_i) \leq rk(K_1) = rk(K)$. Thus the Claim holds.

Suppose the sequence $K_1 \leq K_2 \leq \ldots$ terminates in $m$ steps with $K_m$. By construction this means that $K_m$ is malnormal. The Claim also shows that $K \leq K_m$ is a algebraic extension.

Thus the set $\mathcal{M}(H)$ of algebraic extensions of $K$ which are malnormal in $F(X)$ is nonempty. This sets is also obviously finite. We claim that it has a unique minimal element. Suppose this is not the case and $H_1, H_2$ are two such distinct minimal elements. Since $H_1$ and $H_2$ are malnormal in $F(X)$ and finitely generated, their intersection $H = H_1 \cap H_2$ is also malnormal in $F(X)$, finitely generated and contains $K$. Since $H_1 \neq H_2$ and $H_1, H_2$ are minimal in $\mathcal{M}(H)$, the group $H$ is not an algebraic extension of $K$. That is the extension $K \leq H$ is free.

Let $C \neq 1, K'$ be such that $K \leq K'$, $H = K' * C$ and $C$ is of the largest possible rank. By Lemma 11.9 $K'$ is a algebraic extension of $K$. Since $K'$ is a free factor in $H$, $K'$ is malnormal in $H$ and therefore in $F(X)$. Therefore $K' \in \mathcal{M}(H)$. However, $K' \neq H_1$, $K' \leq H_1$ which contradicts the minimality of $H_1$. $\qquad\square$

**Definition 13.5** (Malnormal closure). Let $K \leq F(X)$ be a finitely generated subgroup. We denote

$$mal(K) := \cap \{ H \mid K \leq H \leq F(X) \text{ and } H \text{ is malnormal in } F(X) \}$$

and refer to $mal(K)$ as the *malnormal closure* of $K$ in $F(X)$.

Note that $F(X)$ is malnormal in $F(X)$, and hence $mal(K)$ is nonempty. Thus $mal(K)$ is the smallest malnormal subgroup of $F(X)$ containing $K$.

**Theorem 13.6.** *Let $K$ be a finitely generated subgroup of $F(X)$.*

*Let $M$ be the unique minimal element among algebraic extensions of $K$ which are malnormal in $F(X)$ (the existence of $H$ follows from Lemma 13.4).*

*Then $M = mal(K)$. In particular $mal(K)$ is finitely generated and the extension $K \leq mal(K)$ is algebraic. Moreover, there is an algorithm which, given a finite set of generators of $K$, outputs a finite set of generators for $mal(K)$.*

*Proof.* Let $K \leq H \leq F(X)$ where $H$ is malnormal in $F(X)$. If $H$ is an algebraic extension of $K$, then $M \leq H$ by the definition of $M$. Suppose now that the extension $K \leq H$ is free. Assume first that $H$ is finitely generated. Let $C \neq 1$ and $K'$ be such that $K \leq K'$, $H = K' * C$ and $C$ is of the largest possible rank. By Lemma 11.9 $K'$ is an algebraic extension of $K$. Since $K'$ is a free factor of $H$, $K$ is malnormal in $H$ and so in $F(X)$. Hence $M \leq K'$ by the definition of $M$ and therefore $M \leq H$.

Assume now that $H$ is not finitely generated, $K \leq H \leq F(X)$ and $H$ is malnormal in $F(X)$. Thus $H$ is a free group of infinite rank. Choose a free basis $Y$ of $H$. Since $K$ is finitely generated, the elements of a finite generating set of $K$ involve only finitely many letters of $Y$. Thus $K$ is contained in a finitely generated free factor $H'$ of $H$. Since $H'$ is a free factor of $H$, it is malnormal in $H$ and so in $F(X)$. Hence $M \leq H'$ by the previous case. Therefore $M \leq H$.

Since $H$ was chosen arbitrarily, this implies that

$$M \leq mal(K)$$

On the other hand $mal(K) \leq M$ by definition of $mal(K)$ since $K \leq M$ and $M$ is malnormal in $F(X)$. Hence $M = mal(K)$, as required.

Thus $K \leq mal(K)$ is an algebraic extension and $mal(K)$ is finitely generated. Moreover, by Theorem 11.3 we can compute all algebraic extensions of $K$ and for each of them check if it is malnormal in $F(X)$ using Corollary 9.11. Since the number of algebraic extensions of $K$ is finite, we can effectively find the unique minimal element among malnormal algebraic extensions of $K$, that is $M$. Thus, we can compute the malnormal closure $mal(K) = M$.

It remains to show that $rk(M) \leq rk(K)$. If $M = K$, there is nothing to prove. Suppose $K \lneqq M$. Hence $K$ is not malnormal in $M$ by the choice of $M$. Therefore there is $g_1 \in M - K$ such that $g_1 K g_1^{-1} \cap K \neq 1$. Arguing exactly as in the proof of Lemma 13.4 we can find a sequence of elements $g_1, g_2, \ldots in M$ such that $K = K_1 \leq K_2 \leq \ldots$ and $K_{i+1} = \langle K_i, g_i \rangle$, where $g_i \in M - K_i$ is such that $g_i K_i g_i^{-1} \cap K_i \neq 1$. Again, by Lemma lem:step each $K_i \leq K_{i+1}$ is an algebraic extension and $rk(K_{i+1}) \leq rk(K_i)$. Hence each $K \leq K_i$ is an algebraic extension and $rk(K_i) \leq rk(K)$. Since the number of algebraic extensions of $K$ in $M$ is finite, any chain of subgroups of this sort terminates. Consider a maximal chain of this type $K = K_1 \leq K_2 \leq \cdots \leq K_s$. Then $K_s$ is malnormal in $M$. Indeed, if not then we can find $g_s \in M - K_s$ with $g_s K_s g_s^{-1} \cap K_s \neq 1$ and extend the chain, contradicting its maximality. Thus $K_s$ is indeed malnormal in $M$ and therefore in $F(X)$. Since $K_s$ is an algebraic extension of $K$, we have $K_s = M$ by the definition of $M$. Therefore $rk(M) \leq rk(K)$, as required.

This completes the proof of Theorem 13.6. □

It turns out that one can use the same argument to analyze *isolators* of finitely generated subgroups of $F(X)$. Recall that a subgroup $H$ of a group $G$ is said to be *isolated* in $G$ if whenever $g^n \in H$, we have $g \in H$. Note that if $G$ is torsion-free and $H$ is not isolated in $G$ then $H$ is not malnormal in $G$. Indeed, suppose $g^n \in H$ but $g \notin H$. Then $1 \neq \langle g^n \rangle = g \langle g^n \rangle g^{-1} \cap \langle g^n \rangle \leq g H g^{-1} \cap H$ and so $H$ is not malnormal. The proof of the following statement is a complete analogue of the proof of Lemma 13.4.

**Lemma 13.7.** *Let $K$ be a finitely generated subgroup of $F(X)$. Then there exists a unique subgroup $H$ of $F(X)$ such that $H$ is minimal among algebraic extensions of $K$ which are isolated in $F(X)$.*

*Proof.* If $K$ is isolated in $F(X)$, the statement is obvious. Assume now that $K$ is not isolated in $F(X)$ (and so not malnormal).

First we observe that there is at least one isolated subgroup which is an algebraic extension of $K$.

Let $K = K_1$ and define a sequence $K_1 \leq K_2 \leq \ldots$ as follows. Since $K_1$ is not isolated, there is $g_2 \in F(X) - K_1$ such that $g_2^{n_2} \in K_1$. Put $K_2 = \langle K_1, g_2 \rangle$. Suppose now that $K_1 \leq K_2 \leq \cdots \leq K_i$ are already constructed. If $K_i$ is isolated in $F(X)$, we terminate the sequence. Otherwise there is $g_{i+1} \in F(X) - K_i$ such that $g_{i+1}^{n_{i+1}} \in K_i$. Put $K_{i+1} = \langle K_i, g_{i+1} \rangle$.

**Claim.** The sequence $K_1 \leq K_2 \leq \ldots$ terminates in a finite number of steps. Moreover, each $K_i$ is a algebraic extension of $K$ and $rk(K_i) \leq rk(K)$.

Note that $g_{i+1} K_i g_{i+1}^{-1} \cap K_i \neq 1$ and so Lemma 13.3 implies that $K_i \leq K_{i+1}$ is a algebraic extension and $rk(K_{i+1}) \leq rk(K_i)$. By transitivity (see Theorem 11.6) this means that $K \leq K_i$ is algebraic for each $i \geq 1$. Since $K_i \neq K_{i+1}$ and there are only finitely many algebraic extensions of $K$, the sequence terminates. Also, obviously $rk(K_i) \leq rk(K_1) = rk(K)$. Thus the Claim holds.

Suppose the sequence $K_1 \leq K_2 \leq \ldots$ terminates in $m$ steps with $K_m$. By construction this means that $K_m$ is isolated in $F(X)$. The Claim also shows that $K \leq K_m$ is an algebraic extension.

Thus the set $\mathcal{I}(H)$ of algebraic extensions of $K$ which are isolated in $F(X)$ is nonempty. This set is also obviously finite. We claim that it has a unique minimal element. Suppose this is not the case and $H_1, H_2$ are two such distinct minimal elements. Since $H_1$ and $H_2$ are isolated in $F(X)$ and finitely generated, their intersection $H = H_1 \cap H_2$ is also isolated in $F(X)$, finitely generated and contains $K$. Since $H_1 \neq H_2$ and $H_1, H_2$ are minimal in $\mathcal{I}(H)$, the group $H$ is not an algebraic extension of $K$. That is the extension $K \leq H$ is free.

Let $C \neq 1, K'$ be such that $K \leq K'$, $H = K' * C$ and $C$ is of the largest possible rank. By Lemma 11.9 $K'$ is a algebraic extension of $K$. Since $K'$ is a free factor in $H$, $K'$ is isolated in $H$ and hence in $F(X)$. Therefore $K' \in \mathcal{I}(H)$. However, $K' \neq H_1$, $K' \leq H_1$ which contradicts the minimality of $H_1$. □

**Definition 13.8** (Isolator). Let $K \leq F(X)$ be a finitely generated subgroup. We define the *isolator* $iso(K)$ of $K$ in $F(X)$ as

$$iso(K) := \cap\{H \mid K \leq H \leq F(X) \text{ and } H \text{ is isolated in } F(X)\}.$$

Thus $iso(K)$ is the smallest isolated subgroup of $F(X)$ containing $K$.

**Theorem 13.9.** *Let $K$ be a finitely generated subgroup of $F(X)$. Let $M$ be the unique minimal element among algebraic extensions of $K$ which are isolated in $F(X)$ (the existence of $H$ follows from Lemma 13.7).*

*Then $iso(K) = M$. In particular, $iso(K)$ is finitely generated, the extension $K \leq iso(K)$ is algebraic and $rk(iso(K)) \leq rk(K)$. Moreover, there exists an algorithm which, given a finite generating set of $K$, produces a free basis of $iso(K)$.*

*Proof.* The proof is exactly the same as that of Theorem 13.6 (for the algorithmic part one needs to use Theorem 13.1 to decide if a particular algebraic extension of $K$ is isolated in $F(X)$). We leave the details to the reader. □

## 14. Ascending chains of subgroups

It is well known that $F(X)$ contains free subgroups of infinite rank if $\#X \geq 2$. Therefore ascending chains of finitely generated subgroups of $F(X)$ do not necessarily terminate. However, by a classical result of M.Takahasi [46] and G.Higman [21] such chains do terminate if the ranks of the subgroups are bounded:

**Theorem 14.1** (Takahasi-Higman). *[46, 21] Let $F = F(X)$ be a free group of finite rank. Let $M \geq 1$ be an integer. Then every strictly ascending chain of subgroups of $F(X)$ of rank at most $M$ terminates.*

*Proof.* Suppose the statement of Theorem 14.1 fails and there is an infinite strictly ascending chain of nontrivial subgroups of $F(X)$

$$(2) \qquad\qquad K_1 \subsetneqq K_2 \subsetneqq \dots$$

where $rk(K_i) \leq m$ for each $i \geq 1$.

Let $\Gamma_i = \Gamma(K_i)$ for $i = 1, 2, \dots$.

We will define the sequence of finite graphs $\Delta_1 \subseteq \Delta_2 \subseteq \dots$ as follows.

**Step 1.** Consider the graphs $\Gamma_{K_2}(K_1), \Gamma_{K_3}(K_1), \dots,$ which are all principal quotients of $\Gamma(K_1)$. Since $\Gamma(K_1)$ is finite, it has only finitely many principal quotients. Therefore for infinitely many values of $i \geq 2$ all $(\Gamma_{K_i}(K_1), 1_{K_i})$ are the same based $X$-digraph which we call $(\Delta_1, v_1)$.

Thus after passing to a subsequence in (2) we assume that $(\Gamma_{K_i}(K_1), 1_{K_i}) = (\Delta_1, v_1)$ for all $i \geq 2$.

**Step $n$.** Suppose the graphs $(\Delta_1, v_1) \subsetneqq (\Delta_2, v_2) \subsetneqq \dots \subsetneqq (\Delta_{n-1}, v_{n-1})$ have already been constructed and the sequence (2) has been modified so that for each $j = 1, \dots, n-1$ we have

$$\Gamma_{K_i}(K_j), 1_{K_i}) = (\Delta_j, v_j) \text{ for all } i \geq j+1.$$

Note that there exists $m \geq n$ such that $\Gamma_{K_m}(K_{n-1}), 1_{K_m}) = (\Delta_{n-1}, v_{n-1})$ is a proper subgraph of $(\Gamma(K_m), 1_{K_m})$. If this is not the case then the sequence of graphs $(\Gamma(K_i), 1_{K_i})$ eventually stabilizes and so $K_i = K_{i+1}$ for some $i$, contrary to our assumptions.

After passing to a subsequence in (2) we may assume that $(\Gamma_{K_n}(K_{n-1}), 1_{K_n}) = (\Delta_{n-1}, v_{n-1})$ is a proper subgraph of $(\Gamma(K_n), 1_{K_n})$.

Since $(\Gamma(K_n), 1_{K_n})$ is a finite graph and has only finitely many principal quotients, there are infinitely many indices $i \geq n+1$ such that all $(\Gamma_{K_i}(K_n), 1_{K_i})$ are the same based $X$-digraph which we call $(\Delta_n, v_n)$. Once again, after passing to a subsequence in (2) we may assume that $(\Gamma_{K_i}(K_n), 1_{K_i}) = (\Delta_n, v_n)$ for all $i \geq n+1$.

Note that $(\Delta_{n-1}, v_{n-1})$ is a subgraph of $(\Delta_n, v_n)$. Indeed, the graph $(\Delta_{n-1}, v_{n-1})$ can be obtained by first mapping $(\Gamma(K_{n-1}), 1_{K_{n-1}})$ into $(\Gamma(K_n), 1_{K_n})$ (under the graph morphism corresponding to $K_{n-1} \leq K_n$) and then mapping the result into $(\Gamma(K_{n+1}), 1_{K_{n+1}})$ (under the graph morphism corresponding to $K_n \leq K_{n+1}$). The full image of the last map from $(\Gamma(K_{n-1}), 1_{K_{n-1}})$ into $(\Gamma(K_{n+1}), 1_{K_{n+1}})$ is precisely $(\Delta_n, v_n)$. Thus $(\Delta_{n-1}, v_{n-1})$ is indeed a subgraph of $(\Delta_n, v_n)$. Moreover it is a proper subgraph. Indeed, suppose $(\Delta_{n-1}, v_{n-1}) = (\Delta_n, v_n)$. Since $(\Gamma_{K_{n+1}}(K_n), 1_{K_{n+1}}) = (\Delta_n, v_n) = (\Delta_{n-1}, v_{n-1})$, we have $K_n \leq L(\Delta_{n-1}, v_{n-1})$. On the other hand by our assumption $(\Delta_{n-1}, v_{n-1})$ is a proper subgraph of $\Gamma(K_n), 1_{K_n})$. Since $\Delta_{n-1}$ is a core graph with respect to $v_{n-1}$, this implies that $L(\Delta_{n-1}, v_{n-1}) \lneqq K_n$, yielding a contradiction. Thus $(\Delta_{n-1}, v_{n-1})$ is a proper subgraph of $(\Delta_n, v_n)$. We now go to the next step.

This procedure gives us an infinite sequence of core graphs $(\Delta_1, v_1) \subsetneqq (\Delta_2, v_2) \subsetneqq \dots$ such that for each $n \geq 1$ $(\Delta_n, v_n)$ is a subgraph of one of the graphs $(\Gamma(K_i), 1_{K_i})$ for some group $K_i$ from the *original* sequence (2). However, since $(\Delta_{n-1}, v_{n-1})$ is a proper subgraph of $(\Delta_n, v_n)$, the group $L(\Delta_{n-1}, v_{n-1})$ is a proper free factor of the group $L(\Delta_n, v_n)$ for every $n \geq 1$. Therefore $rk(L(\Delta_n, v_n)) \geq n$ for each $n \geq 1$. However, $(\Delta_{M+1}, v_{M+1})$ is a subgraph of some $(\Gamma(K_i), 1_{K_i})$ and therefore $L(\Delta_{M+1}, v_{M+1})$ is a free factor of $K_i$. Since the rank of $L(\Delta_{M+1}, v_{M+1})$ is at least $M + 1$, this implies that the rank of $K_i$ is at least $M + 1$ as well. This contradicts our assumption in Theorem 14.1 that the ranks of all subgroups in (2) are at most $M$.  □

## 15. Descending chains of subgroups

In this section we will show that graph-theoretic methods developed earlier allow one to work with descending as well as ascending chains of subgroups. Namely, we will prove two classical results due to M.Takahasi [46].

**Theorem 15.1.** [46] *Suppose $K_i$ are subgroups of $F(X)$ which form an infinite strictly descending chain*

$$K_1 \gneqq K_2 \gneqq \dots$$

*Put $K_\infty = \cap_{n=1}^\infty K_n$. Let $K$ be a finitely generated free factor of $K_\infty$. Then $K$ is a free factor in all but finitely many $K_i$.*

*Proof.* Note that $\Gamma(K)$ is a finite graph and therefore it has only finitely many principal quotients.

Therefore there is an infinite strictly increasing sequence of indices $(i_n)_{n=1}^\infty$ such that for all $n \geq 1$ the based graphs $(\Gamma_{K_{i_n}}(K), 1_{K_{i_n}})$ are the same graph $(\Delta, v)$. Let $H = L(\Delta, v)$ so that $(\Gamma(H), 1_H) = (\Delta, v)$.

We claim that in fact $(\Gamma_{K_j}(K), 1_{K_j}) = (\Delta, v)$ for each $j \geq i_1$. Indeed, suppose $i_n < j < i_{n+1}$. Since $i_n < j$ and $K_j \leq K_{i_n}$ the graph $(\Gamma_{K_j}(K), 1_{K_j})$ is an epimorphic image of the graph $(\Gamma_{K_{i_n}}(K), 1_{K_{i_n}}) = (\Delta, v)$. Therefore $L(\Gamma_{K_j}(K), 1_{K_j}) \leq L(\Delta, v)$. On the other hand $j < i_{n+1}$, $K_{i_{n+1}} \leq K_j$ and therefore the graph $(\Delta, v) = (\Gamma_{K_{i_{n+1}}}(K), 1_{K_{i_{n+1}}})$ is an epimorphic image of the graph $(\Gamma_{K_j}(K), 1_{K_j})$. Hence $L(\Delta, v) \leq L(\Gamma_{K_j}(K), 1_{K_j})$. Thus $L(\Gamma(H), 1_H) = L(\Delta, v) = L(\Gamma_{K_j}(K))$. Since $\Delta$ is a finite core graph with respect to $v$, this implies that $(\Gamma(H), 1_H) = (\Delta, v) = (\Gamma_{K_j}(K), 1_{K_j})$, as claimed.

Since $(\Gamma(H), 1_H) = (\Delta, v)$ is a subgraph of all $(\Gamma(K_j), 1_{K_j})$ for $j \geq i_1$, the group $H$ is a free factor in all $K_j$ for $j \geq i_1$. Therefore $H \leq \cap_{i=1}^\infty K_i = K_\infty$.

Also, since $(\Delta, v)$ is an epimorphic image of $(\Gamma(K), 1_K)$, we have $K = L(\Gamma(K), 1_K) \leq L(\Delta, v) = H$.

Recall that $K$ is a finitely generated free factor of $K_\infty$ and $K \leq H \leq K_\infty$. Therefore by Lemma 10.7 $K$ is a free factor of $H$. Since $H$ is a free factor of all but finitely many $K_i$, the same holds for $K$.

This completes the proof of Theorem 15.1 is proved.  □

**Theorem 15.2.** [46] *Let $F(X)$ be a free group of finite rank. Let*

$$K_1 \gneqq K_2 \gneqq \dots$$

*be an infinite strictly decreasing chain of subgroups in $F(X)$. Suppose also that $rk(K_i) \leq M$ for each $i \geq 1$. Put $K_\infty = \cap_{i=1}^\infty K_i$.*

*Then:*
*(a) $K_\infty$ is a free factor in all but finitely many $K_i$.*
*(b) The group $K_\infty$ is finitely generated and its rank is at most $M - 1$.*

*Proof.* Note that (a) implies (b). In fact suppose (a) holds. Since $K_\infty$ is a free factor in some $K_i$, the free group $K_\infty$ is finitely generated and has rank at most $M$. Suppose that in fact $K_\infty$ is free of rank $M$. Since a proper free factor of a free group always has small rank and $rk(K_i) \leq M$ this implies that $K_\infty = K_i$. However $K_{i+1} \subsetneqq K_i$ and therefore $K_{i+1}$ is not contained in $K_\infty$. This contradicts the fact that $K_\infty = \cap_{n=1}^\infty K_n$. Thus $rk(K_\infty) \le M - 1$ as required.

We will now establish part (a) of the theorem.

Suppose first that $K_\infty$ is finitely generated. Since $K_\infty$ is a free factor of itself, Theorem 15.1 implies that $K_\infty$ is a free factor in all but finitely many $K_i$, as required.

Suppose now that $K_\infty$ is a free group of infinite rank. Then there are finitely generated free factors of $K_\infty$ of arbitrarily large rank. This contradicts the conclusion of Theorem 15.1 since the ranks of $K_i$ are bounded. $\qquad\square$

Department of Mathematics, University of Illinois at Urbana-Champaign, 1409 West Green Street, Urbana, IL 61801, USA
*E-mail address*: `kapovich@math.uiuc.edu`

Department of Mathematics, City College of CUNY, New York, NY 10031, USA
*E-mail address*: `alexei@rio.sci.ccny.cuny.edu`